\documentclass[11pt]{article}
\usepackage{amsmath, upgreek}
\usepackage{amssymb}
\usepackage{color}
\usepackage{amscd}
\usepackage{xspace}
\usepackage{verbatim}
\usepackage{graphicx}
\newcommand{\mysection}[1]{
\section{#1}\setcounter{equation}{0}}
\title{\bf Estimates of solutions of elliptic equations
 with a source reaction term involving the product of the function and its gradient}%%
\author{{\bf Marie-Fran\c{c}oise Bidaut-V\'eron\footnote{\noindent Laboratoire de Math\'{e}matiques et Physique Th\'{e}orique,
Universit\'e de Tours, 37200 Tours, France. E-mail: veronmf@univ-tours.fr},} \\{\bf Marta Garcia-Huidobro \footnote{\noindent
Departamento de Matematicas, Pontifica Universidad Catolica de Chile
Casilla 307, Correo 2, Santiago de Chile. E-mail: mgarcia@mat.puc.cl}}\\
 {\bf Laurent V\'eron \footnote{\noindent
Laboratoire de Math\'{e}matiques et Physique Th\'{e}orique, Universit\'e de Tours, 37200 Tours, France. E-mail: veronl@univ-tours.fr}}\\[2mm]
}%%��

\date{}
\begin{document}
 \maketitle
% \noindent{\small {\bf Abstract} We study the existence and uniqueness of  solutions of $\partial_tu-\Delta u+u^q=0$ ($q>1$) in $\Omega\times (0,\infty)$ where $\Omega\subset\mathbb R^N$ is a domain with a compact boundary, subject to the conditions $u=f\geq 0$ on $\partial\Omega\times (0,\infty)$ and the initial condition $\lim_{t\to 0}u(x,t)=\infty$. By means of Brezis' theory of maximal monotone operators in Hilbert spaces, we construct a minimal solution when $f=0$, whatever is the regularity of the boundary of the domain. When $\partial\Omega$ satisfies the parabolic Wiener criterion and $f$ is continuous, we construct a maximal solution and prove that it is the unique solution which blows-up at $t=0$.
% }

% \noindent
% {\it \footnotesize 1991 Mathematics Subject Classification}. {\scriptsize
% 35K60}.\\
% {\it \footnotesize Key words}. {\scriptsize Parabolic equations, singular solutions, semi-groups of contractions, maximal monotone operators, Wiener criterion.}
% \vspace{1mm}
% \hspace{.05in}

%% FONT commands
\newcommand{\txt}[1]{\;\text{ #1 }\;}%% Used in math only
\newcommand{\tbf}{\textbf}%% Bold face. Usage: \tbf{...}
\newcommand{\tit}{\textit}%% Italic
\newcommand{\tsc}{\textsc}%% Small caps
\newcommand{\trm}{\textrm}
\newcommand{\mbf}{\mathbf}%% Math bold
\newcommand{\mrm}{\mathrm}%% Math Roman
\newcommand{\bsym}{\boldsymbol}%% Bold math symbol
%%Macros for changing font size in math.
\newcommand{\scs}{\scriptstyle}%% as in subscript
\newcommand{\sss}{\scriptscriptstyle}%% as in sub-subscript
\newcommand{\txts}{\textstyle}
\newcommand{\dsps}{\displaystyle}
%%Macros for changing font size in text.
\newcommand{\fnz}{\footnotesize}
\newcommand{\scz}{\scriptsize}
%%\tiny<\scz<\fsz<\small<\large<\Large<\huge<\Huge
%%%%%%%%%%%%
%%%%%%%%%%%%
%% EQUATION commands
\newcommand{\be}{\begin{equation}}
\newcommand{\bel}[1]{\begin{equation}\label{#1}}
\newcommand{\ee}{\end{equation}}%% This macro does not work with amstex.
\newcommand{\eqnl}[2]{\begin{equation}\label{#1}{#2}\end{equation}}%%use not advisable; confusing
\newcommand{\barr}{\begin{eqnarray}}
\newcommand{\earr}{\end{eqnarray}}
\newcommand{\bars}{\begin{eqnarray*}}
\newcommand{\ears}{\end{eqnarray*}}
\newcommand{\nnu}{\nonumber \\}
%%%%%%%%%%%%%%%
%% Unnumbered THEOREM env.
%% New env. to be used for unnumbered theorem, lemma etc. (but with specified name)
\newtheorem{subn}{\name}
\renewcommand{\thesubn}{}
\newcommand{\bsn}[1]{\def\name{#1}\begin{subn}}
\newcommand{\esn}{\end{subn}}
%%%%%%%%%%%%%%
%% NUMBERED THEOREM env.
%% Environments: theorem, lemma, corollary defintion and related commands,
%% designed to provide consecutive numbering of these forms.
\newtheorem{sub}{\name}[section]
\newcommand{\dn}[1]{\def\name{#1}}   %used in conjuction with sub or subn.
\newcommand{\bs}{\begin{sub}}
\newcommand{\es}{\end{sub}}
\newcommand{\bsl}[1]{\begin{sub}\label{#1}}
%% the above must be preceeded by \dn (name definition),
%% however this is superceded by the list of commands bth etc.  below.
%%%%%%%%%%%%
%% NUMBERED THEOREM env. (cont.)
%% List of commands derived from 'sub' env. for theorem, lemma etc.
%% designed to provide consecutive numbering of these forms.
\newcommand{\bth}[1]{\def\name{Theorem}
\begin{sub}\label{t:#1}}
\newcommand{\blemma}[1]{\def\name{Lemma}
\begin{sub}\label{l:#1}}
\newcommand{\bcor}[1]{\def\name{Corollary}
\begin{sub}\label{c:#1}}
\newcommand{\bdef}[1]{\def\name{Definition}
\begin{sub}\label{d:#1}}
\newcommand{\bprop}[1]{\def\name{Proposition}
\begin{sub}\label{p:#1}}
%%%%%%%%%%%%%%%%%%%%%%%%%%%%%%%%%%
%% RERERENCE commands.
%% \newcommand{\R}[1]{(\ref{#1})}
\newcommand{\R}{\eqref}
\newcommand{\rth}[1]{Theorem~\ref{t:#1}}
\newcommand{\rlemma}[1]{Lemma~\ref{l:#1}}
\newcommand{\rcor}[1]{Corollary~\ref{c:#1}}
\newcommand{\rdef}[1]{Definition~\ref{d:#1}}
\newcommand{\rprop}[1]{Proposition~\ref{p:#1}}
%%%%%%%%%%%
%% ARRAY commands.
\newcommand{\BA}{\begin{array}}
\newcommand{\EA}{\end{array}}
\newcommand{\BAN}{\renewcommand{\arraystretch}{1.2}
\setlength{\arraycolsep}{2pt}\begin{array}}
\newcommand{\BAV}[2]{\renewcommand{\arraystretch}{#1}
\setlength{\arraycolsep}{#2}\begin{array}}
%Note: The first variable gives the amount of stretching: (#1) x default.
%For instance #1=1.2 means a 20% stretching. The second variable should be
%written for instance in the form  4pt ; here the default is 5pt
%\newcommand{\EAN}{\end{array}\setlength{\arraycolsep}{5pt}}
\newcommand{\BSA}{\begin{subarray}}
\newcommand{\ESA}{\end{subarray}}
%Note: These are used in subscripts as well as superscripts. They work essentially
%% like 'array'.
\newcommand{\BAL}{\begin{aligned}}
\newcommand{\EAL}{\end{aligned}}
\newcommand{\BALG}{\begin{alignat}}
\newcommand{\EALG}{\end{alignat}}%% the abbrev. does not work with latex2e
\newcommand{\BALGN}{\begin{alignat*}}
\newcommand{\EALGN}{\end{alignat*}}%% the abbrev. does not work with latex2e
%% The 'aligned' environment must be placed inside an 'equation' env.
%% in the same way as the array.
%% One could use also the 'align' env. or the 'alignat' env.
%% However in this case each line is numbered, unless '\notag' is used.
%% The 'alignat'
%% has a slightly different format (the number of columns must be specified in advance)
%% but it has the advantage that the distance between columns is at our disposition.
%% (The default would be zero distance.) Using 'alignat*' we can have the advantages
%% of alignat plus the situation where separate lines are not numbered.
%% However in this case there is no numbering at all (unless we provide a tag).
%%%%%%%%%%
%% PROOF, REMARK etc.
\newcommand{\note}[1]{\textit{#1.}\hspace{2mm}}
\newcommand{\Proof}{\note{Proof}}
\newcommand{\qeda}{\hspace{10mm}\hfill $\square$}
\newcommand{\qed}{\\
${}$ \hfill $\square$}
\newcommand{\Remark}{\note{Remark}}
%%%%%%%% Style command.
\newcommand{\modin}{$\,$\\[-4mm] \indent}
%% To be used after \mysection in order to start new line with \indent.
%%%%%%%%%%%%
%% MATHEMATICAL symbols
\newcommand{\forevery}{\quad \forall}
\newcommand{\set}[1]{\{#1\}}
\newcommand{\setdef}[2]{\{\,#1:\,#2\,\}}
\newcommand{\setm}[2]{\{\,#1\mid #2\,\}}
%% Arrows
\newcommand{\mt}{\mapsto}
\newcommand{\lra}{\longrightarrow}
\newcommand{\lla}{\longleftarrow}
\newcommand{\llra}{\longleftrightarrow}
\newcommand{\Lra}{\Longrightarrow}
\newcommand{\Lla}{\Longleftarrow}
\newcommand{\Llra}{\Longleftrightarrow}
\newcommand{\warrow}{\rightharpoonup}
%% Brackets, delimiters
\newcommand{
\paran}[1]{\left (#1 \right )}%% adjustable parantheses
\newcommand{\sqbr}[1]{\left [#1 \right ]}%% adjustable square brackets
\newcommand{\curlybr}[1]{\left \{#1 \right \}}%% adjustable curly brackets
\newcommand{\abs}[1]{\left |#1\right |}%% adjustable vertical delimiters
\newcommand{\norm}[1]{\left \|#1\right \|}%% adjustable norm
\newcommand{
\paranb}[1]{\big (#1 \big )}%% non-adjustable parantheses (big)
\newcommand{\lsqbrb}[1]{\big [#1 \big ]}%% non-adjustable square brackets (big)
\newcommand{\lcurlybrb}[1]{\big \{#1 \big \}}%% non-adjustable curly brackets (big)
\newcommand{\absb}[1]{\big |#1\big |}%% non-adjustable vertical delimiters (big)
\newcommand{\normb}[1]{\big \|#1\big \|}%% non-adjustable norm (big)
\newcommand{
\paranB}[1]{\Big (#1 \Big )}%% non-adjustable parantheses (Big)
\newcommand{\absB}[1]{\Big |#1\Big |}%% non-adjustable vertical delimiters (Big)
\newcommand{\normB}[1]{\Big \|#1\Big \|}%% non-adjustable norm (Big)
\newcommand{\produal}[1]{\langle #1 \rangle}%% the pairing of X' and X

%%%%%%%%%%%%%%%%%
%% Adjustable parantheses etc. in a different DEFINITION format.
%\def\adp(#1){\left (#1 \right )}%% adjustable parantheses
%\def\adsb(#1){\left [#1\right ]}%% adjustable square brackets
%\def\adcb(#1){\left \{#1\right \}}%% adjustable curly brackets
%\def\abs|#1|{\left |#1\right |}%% adjustable vertical delimiters
%%%%%%%%%%%%%%%%
%% More mathematical symbols
\newcommand{\thkl}{\rule[-.5mm]{.3mm}{3mm}}
\newcommand{\thknorm}[1]{\thkl #1 \thkl\,}
\newcommand{\trinorm}[1]{|\!|\!| #1 |\!|\!|\,}
\newcommand{\bang}[1]{\langle #1 \rangle}%% angle bracket
\def\angb<#1>{\langle #1 \rangle}%% angle bracket
%% The two last lines yield the same result.
%% The second is used as follows: \angb<a,b>
\newcommand{\vstrut}[1]{\rule{0mm}{#1}}
\newcommand{\rec}[1]{\frac{1}{#1}}
%% OPERATOR names.
%% OPERATOR names.
\newcommand{\opname}[1]{\mbox{\rm #1}\,}
\newcommand{\supp}{\opname{supp}}
\newcommand{\dist}{\opname{dist}}
\newcommand{\myfrac}[2]{{\displaystyle \frac{#1}{#2} }}
\newcommand{\myint}[2]{{\displaystyle \int_{#1}^{#2}}}
\newcommand{\mysum}[2]{{\displaystyle \sum_{#1}^{#2}}}
\newcommand {\dint}{{\displaystyle \myint\!\!\myint}}%%%%%%%%%%
%%%%%%% SPACE commands
\newcommand{\q}{\quad}
\newcommand{\qq}{\qquad}
\newcommand{\hsp}[1]{\hspace{#1mm}}
\newcommand{\vsp}[1]{\vspace{#1mm}}
%%%%%%%%%%%
%% ABREVIATIONS
\newcommand{\ity}{\infty}
\newcommand{\prt}{\partial}
\newcommand{\sms}{\setminus}
\newcommand{\ems}{\emptyset}
\newcommand{\ti}{\times}
\newcommand{\pr}{^\prime}
\newcommand{\ppr}{^{\prime\prime}}
\newcommand{\tl}{\tilde}
\newcommand{\sbs}{\subset}
\newcommand{\sbeq}{\subseteq}
\newcommand{\nind}{\noindent}
\newcommand{\ind}{\indent}
\newcommand{\ovl}{\overline}
\newcommand{\unl}{\underline}
\newcommand{\nin}{\not\in}
\newcommand{\pfrac}[2]{\genfrac{(}{)}{}{}{#1}{#2}}% frac with parantheses.
%%%%%%%%%%%
%%%%%%%%%%%%%

%%Macros for Greek letters.
\def\ga{\alpha}     \def\gb{\beta}       \def\gg{\gamma}
\def\gc{\chi}       \def\gd{\delta}      \def\ge{\epsilon}
\def\gth{\theta}                         \def\vge{\varepsilon}
\def\gf{\phi}       \def\vgf{\varphi}    \def\gh{\eta}
\def\gi{\iota}      \def\gk{\kappa}      \def\gl{\lambda}
\def\gm{\mu}        \def\gn{\nu}         \def\gp{\pi}
\def\vgp{\varpi}    \def\gr{\rho}        \def\vgr{\varrho}
\def\gs{\sigma}     \def\vgs{\varsigma}  \def\gt{\tau}
\def\gu{\upsilon}   \def\gv{\vartheta}   \def\gw{\omega}
\def\gx{\xi}        \def\gy{\psi}        \def\gz{\zeta}
\def\Gg{\Gamma}     \def\Gd{\Delta}      \def\Gf{\Phi}
\def\Gth{\Theta}
\def\Gl{\Lambda}    \def\Gs{\Sigma}      \def\Gp{\Pi}
\def\Gw{\Omega}     \def\Gx{\Xi}         \def\Gy{\Psi}

%%Macros for calligraphic letters.
\def\CS{{\mathcal S}}   \def\CM{{\mathcal M}}   \def\CN{{\mathcal N}}
\def\CR{{\mathcal R}}   \def\CO{{\mathcal O}}   \def\CP{{\mathcal P}}
\def\CA{{\mathcal A}}   \def\CB{{\mathcal B}}   \def\CC{{\mathcal C}}
\def\CD{{\mathcal D}}   \def\CE{{\mathcal E}}   \def\CF{{\mathcal F}}
\def\CG{{\mathcal G}}   \def\CH{{\mathcal H}}   \def\CI{{\mathcal I}}
\def\CJ{{\mathcal J}}   \def\CK{{\mathcal K}}   \def\CL{{\mathcal L}}
\def\CT{{\mathcal T}}   \def\CU{{\mathcal U}}   \def\CV{{\mathcal V}}
\def\CZ{{\mathcal Z}}   \def\CX{{\mathcal X}}   \def\CY{{\mathcal Y}}
\def\CW{{\mathcal W}} \def\CQ{{\mathcal Q}}
%%%%%
%%Macros for 'blackboard' letters (See (27) for display.)
\def\BBA {\mathbb A}   \def\BBb {\mathbb B}    \def\BBC {\mathbb C}
\def\BBD {\mathbb D}   \def\BBE {\mathbb E}    \def\BBF {\mathbb F}
\def\BBG {\mathbb G}   \def\BBH {\mathbb H}    \def\BBI {\mathbb I}
\def\BBJ {\mathbb J}   \def\BBK {\mathbb K}    \def\BBL {\mathbb L}
\def\BBM {\mathbb M}   \def\BBN {\mathbb N}    \def\BBO {\mathbb O}
\def\BBP {\mathbb P}   \def\BBR {\mathbb R}    \def\BBS {\mathbb S}
\def\BBT {\mathbb T}   \def\BBU {\mathbb U}    \def\BBV {\mathbb V}
\def\BBW {\mathbb W}   \def\BBX {\mathbb X}    \def\BBY {\mathbb Y}
\def\BBZ {\mathbb Z}

%%Macros for Ghotic (Fraktur) letters.
\def\GTA {\mathfrak A}   \def\GTB {\mathfrak B}    \def\GTC {\mathfrak C}
\def\GTD {\mathfrak D}   \def\GTE {\mathfrak E}    \def\GTF {\mathfrak F}
\def\GTG {\mathfrak G}   \def\GTH {\mathfrak H}    \def\GTI {\mathfrak I}
\def\GTJ {\mathfrak J}   \def\GTK {\mathfrak K}    \def\GTL {\mathfrak L}
\def\GTM {\mathfrak M}   \def\GTN {\mathfrak N}    \def\GTO {\mathfrak O}
\def\GTP {\mathfrak P}   \def\GTR {\mathfrak R}    \def\GTS {\mathfrak S}
\def\GTT {\mathfrak T}   \def\GTU {\mathfrak U}    \def\GTV {\mathfrak V}
\def\GTW {\mathfrak W}   \def\GTX {\mathfrak X}    \def\GTY {\mathfrak Y}
\def\GTZ {\mathfrak Z}   \def\GTQ {\mathfrak Q}

\font\Sym= msam10 % special symbols
\def\SYM#1{\hbox{\Sym #1}}
\newcommand{\bdw}{\prt\Gw\xspace}
\date{}
\maketitle\medskip

\noindent{\small {\bf Abstract} We study local and global properties of positive solutions of $-\Gd u=u^p\abs{\nabla u}^q$ in a domain $\Gw$ of $\BBR^N$, in the range $p+q>1$, $p\geq 0$, $0\leq q< 2$. We first prove  local Harnack inequality and nonexistence of positive solutions in $\BBR^N$ when
$p(N-2)+q(N-1)<N$. Using a direct Bernstein method we obtain a first range of values of $p$ and $q$ in which
$u(x)\leq c(\dist(x,\prt\Gw))^{\frac{q-2}{p+q-1}}$ This holds in particular if $p+q<1+\frac{4}{N-1}$. Using an integral Bernstein method we obtain a wider range of values of $p$ and $q$ in which all the global solutions are constants. Our result contains Gidas and Spruck nonexistence result as a particular case. We also study solutions under the form $u(x)=r^{\frac{q-2}{p+q-1}}\gw(\gs)$. We prove existence, nonexistence and rigidity of the spherical component $\gw$ in some range of values of $N$, $p$ and $q$.
}\smallskip

\noindent
{\it \footnotesize 2010 Mathematics Subject Classification}. {\scriptsize 35J62, 35B08, 68Ð04}.\\
{\it \footnotesize Key words}. {\scriptsize elliptic equations; Bernstein methods; gradient estimates; global solutions; bifurcations.
}
\tableofcontents
\vspace{1mm}
\hspace{.05in}
%\abstract
\medskip
%%%%%%%%%%%%%%%%%%%%%%%%%%%%%%%%%%%%%%%%%%%%%%%%%%%%%%%%%%%%%%%%%%%%%%%%%%%%%%%%%%%%%%%%%%%%%%%%%%%%%%%%%%%%%%%%%%%%%%%%%%%%%%%%%%%%%%%%%%%%%%%%%%%%%%%%%%%%%%%%%%%%%%%%%%%%%%%%%%%%%%%%%%%%%
%%%%%%%%%%%%%%%%%%%%%%%%%%%%%%%%%%%%%%%%%%%%%%%%%%%%%%%%%%%%%%%%%%%%%%%%%%%%%%%IN ADOMAIN%%%%%%%%%%%%%%%%%%%%%%%%%%%%%%%%%%%%%%%%%%%%%%%%%%%%%%%%%%%%%%%%%%%%%%%%%%%%%%%%%%%%%%
%%%%%%%%%%%%%%%%%%%%%%%%%%%%%%%%%%%%%%%%%%%%%%%%%%%%%%%%%%%%%%%%%%%%%%%%%%%%%%%%%%%%%%%%%%%%%%%%%%%%%%%%%%%%%%%%%%%%%%%%%%%%%%
%%%%%%%%%%%%%%%%%%%%%%%%%%%%%%%%%%%%%%%%%%%%%%%%%%%%%%%%%%%%%%%%%%%%%%%%%%%%%%%%%%%%%%%%%%%%%%%%%%%%%%%%%%%%%%%%%%%%%%%%%%%%%%
%%%%%%%%%%%%%%%%%%%%%%%%%%%%%%%%%%%%%%%%%%%%%%%%%%%%%%%%%%%%%%%%%%%%%%%%%%%%%%%%%%%%%%%%%%%%%%%%%%%%%%%%%%%%%%%%%%%%%%%%%%%%%%%%%%%%%%%%%%%%%%%%%%%%%%%%%%%%%%%%%%%%%%%%%%%%%%%%%%%%
%%%%%%%%%%%%%%%%%%%%%%%%%%%%%%%%%%%%%%%%%%%%%%%%%%%%%%%%%%%%%SECTION--INTRODUCTION%%%%%%%%%%%%%%%%%%%%%%%%%%%%%%%%%%%%%%%%%%%%%%%%%%%%%%%%%%%%%%%%%%%%%%%%%%%%%%%%%%%%%%%%%%%%%%%%%%%%%%%%%%%%%%%%%%%%%%%%%%%%%%%%%%%%%%%%%%%%%%%%%%%%%%%%%%%%%%%%%%%%
\mysection{Introduction}
The aim of this article is  to study local and global properties of positive solutions of the following type of equations
\bel{I-0}
-\Gd u=u^p\abs{\nabla u}^q,
\ee
in $\Gw\setminus\{0\}$ where $\Gw$ is an open subset of $\BBR^N$ containing $0$, $p$ and $q$ are real exponents. In many cases we will assume  the superlinearity of the right-hand side, i.e. $p+q-1>0$ and $0\leq q \leq 2$.
Equation $(\ref{I-0})$ is invariant under the action of the transformations $T_\gs$ defined for $\gs>0$ by
\bel{I-3}
T_\gs[u](x)=\gs^{\frac{2-q}{p+q-1}}u(\gs x).
\ee
If we look for radial positive solutions under the form $u(x)=\Gl\abs x^{-\gg}$ we find, if $q<2$ and $p+q-1>0$, $\gg:=\gg_{p,q}=\frac{2-q}{p+q-1}$ and
\bel{I-4}
\Gl:=\Gl_{N,p,q}=\gg_{p,q}^{\frac{1-q}{p+q-1}}\left(N-\myfrac{2p+q}{p+q-1}\right)^{\frac{1}{p+q-1}}.
\ee
However this last quantity  exists if and only if the exponents belong to the {\it supercritical range}, that is when
\bel{I-5}
(N-2)p+(N-1)q>N.
\ee
In the {\it subcritical range} of exponents i.e. when
\bel{I-6}
(N-2)p+(N-1)q<N,
\ee
we prove that Serrin's classical results  (see \cite{Se1}, \cite{Se2}) can be applied. We obtain a local Harnack inequality and an {\it a priori} estimate for positive solution
$u$ in $B_R\setminus\{0\}$ under the form
\bel{I-7}
u(x)+\abs x\abs{\nabla u(x)}\leq c\abs x^{2-N}\qquad\forall x\text{ s.t. }0<\abs x\leq \frac{R}{2},
\ee
 with a constant $c$ depending on $u$.  \medskip

 \nind{\bf Theorem A }{\it  Let $\Gw\subset\BBR^N$ be a domain containing $0$, $N\geq 3$, $p\geq 0$, $0\leq q\leq 2$ and assume $(\ref{I-6})$ holds. If $u\in C^2(\Gw\setminus\{0\})$ is a positive solution of $(\ref{I-0})$ in $\Gw\setminus\{0\}$, then estimate $(\ref{I-7})$ holds in a neighborhood of $0$}.\smallskip

If $(\ref{I-6})$ holds all nonnegative functions $u\in C^1(\BBR^N)$ satisfying
\bel{I-7-1}
-\Gd u\geq u^p|\nabla u|^q\geq 0\qquad\text{in }\;\BBR^N,
\ee
are constant (see \cite{BGMQ} for related results). This result is due to Mitidieri and Pohozaev \cite{MiPo}. For the sake of completeness we give a slightly different proof which introduces the techniques we developed throughout our article.

 \smallskip

Our main results deal with the supercritical range. We prove {\it a priori} estimates of positive solutions of $(\ref{I-0})$ in a punctured domain and existence of ground states in $\BBR^N$. There are two approaches for obtaining these results. The {\it direct Bernstein method} and the {\it integral Bernstein method} popularized by Lions \cite{Lio} and Gidas and Spruck in \cite{GS} respectively.
Both methods are based upon differentiating the equation. The direct Bernstein method relies on obtaining pointwise estimates of the gradient {\it without any integration}.  Thanks to highly nontrivial algebraic manipulations with an intensive use of Young's inequality, we  prove that the norm of the gradient of a power of a solution satisfies an  elliptic inequality with a superlinear absorption term, which allows to use the Keller-Osserman comparison  method. Our main result in this framework is the following:\medskip

\nind{\bf Theorem B }{\it  Let   $N\geq 2$, $0\leq q< 2$ and $p\geq 0$ be such that $p+q-1>0$. If $u$ is a positive solution of $(\ref{I-0})$ in $B_{R}$ and one of the following assumptions is fulfilled,\smallskip

\nind (i) $1\leq p< \frac{N+3}{N-1}$ and $p+q-1<\frac{4}{N-1}$,\smallskip

\nind (ii) $0\leq p<1$ and $ p+q-1<\frac{(p+1)^2}{p(N-1)}$.

\smallskip

\nind Then there exist positive constants $a=a(N,p,q)$ and $c_1=c_1(N,p,q)$ such that
\bel{I-8}
\abs{\nabla u^a(0)}\leq c_1R^{-1-a\frac{2-q}{p+q-1}}.
\ee
}\smallskip

Notice that there always holds $\frac{4}{N-1}\leq \frac{(p+1)^2}{p(N-1)}$. The value of the exponent $a$ is not easy to compute, however, in several applications this difficulty can be bypassed. As a consequence of $(\ref{I-8})$ there holds,\medskip

\nind{\bf Corollary B-1 }{\it Under the assumptions on $N$, $p$ and $q$ of Theorem B, any positive solution of $(\ref{I-0})$ in $\BBR^N$ is constant.
}

\medskip

Another consequence is the following,\medskip

\nind{\bf Corollary B-2 }{\it  Let $\Gw$ be a smooth domain in $\BBR^N$, $N\geq 2$, with a bounded boundary, $0\leq q< 2$ and $p\geq 0$ such that $p+q-1>0$ and assume one of the assumptions (i)-(ii) of Theorem B holds. If $u$ is a positive solution of $(\ref{I-0})$ in $\Gw$ there exists $d_0>0$ depending on $\Gw$ and $c_2=c_2(N,p,q)>0$ such that
\bel{I-8-1}
u(x)\leq c_2\left(\left(\dist (x,\prt\Gw)\right)^{-\frac{2-q}{p+q-1}}+\max\{ u(z):\dist (z,\Gw)=d_0\}\right).
\ee
}

The integral Bernstein method has been initiateded in \cite{GS} in order to prove that positive solutions of Lane-Emden equation, i.e. equation $(\ref{I-0})$ with $q=0$, satisfies Harnack inequality near a singularity. Therein they introduced the equation satisfied by a the norm of the gradient of a power of a solution and proved that the term $u^{p-1}$ in this equation verifies some local integral estimate which allows to use Serrin's results \cite{Se1} for getting this inequality. Our method here is to start from the equation satisfied by the norm of the gradient of a power of a solution of $(\ref{I-0})$, to multiply it by a suitable power of this norm and to obtain estimates of the $L^r$-norm of the gradient of the solutions in balls for $r$ large enough. These integral estimates allow us to prove the non-existence of non-constant global solutions.\medskip

 \nind{\bf Theorem C }{\it Assume $p\geq 0$, $0\leq q<2$ and define the polynomial $G$ by
  \bel{I-B-1}\BA {lll}
   G(p,q)= \left((N-1)^2q + N-2\right)p^2 + b(q)p-Nq^2, \\[2mm]
\text {where} \quad b(q)=N(N-1)q^2-(N^2 + N-1)q-N-2.
\EA\ee
 If the couple $(p,q)$ satisfies the inequality $G(p,q)<0$, then all the positive solutions of $(\ref{I-0})$ in $\BBR^N$ are constant.
  } \medskip

In the range of $p$ and $q$, the condition $G_N(p,q)<0$ is equivalent to
  \bel{I-B-2}\BA {lll}\displaystyle
 0\leq p<p_c(q):=\myfrac{-b(q)+\sqrt{b^2(q)+4Nq^2\left((N-1)^2q+N-2\right)}}{2\left((N-1)^2q+N-2\right)}.
\EA\ee
It is noticeable that the minimum of $p$ on the curve $G(p,q)=0$, $0<q<2$ is smaller that $1$ whenever $N\geq 9$. If $q=0$, the above reads reads
  \bel{I-B-3}\BA {lll}\displaystyle
 0\leq p<p_c(0):=\myfrac{N+2}{N-2},
\EA\ee
 which is the well known condition obtained by Gidas and Spruck in \cite{GS}. Furthermore, it can be verified that the domain of $(p,q)$ in which Theorem B applies is included into the set of $(p,q)$ where $G(p,q)<0$. Our proof is extremely technical and necessitates a long appendix in which many algebraic computations are carried out.\smallskip

 If we just look for radial solutions we obtain an optimal result, namely:\smallskip

 \nind{\bf Theorem D }{\it There exist non-constant radial positive solutions of $(\ref{I-0})$ in $\BBR^N$ if and only if $p\geq 0$, $0\leq q<1$ and
    \bel{I-B-3'} p(N-2)+q(N-1)\geq N+\myfrac{2-q}{1-q}.
 \ee
If equality holds in $(\ref{I-B-3'})$, there exists an explicit one parameter family of positive radial solutions of $(\ref{I-0})$ in $\BBR^N$
under the form
    \bel{I-B-3''}
    u_c(r)=c\left(Kc^{\frac{(2-q)^2}{(N-2)(1-q)}}+r^{\frac{2-q}{1-q}}\right)^{-\frac{(N-2)(1-q)}{2-q}},
 \ee
for any $c>0$ and some $K=K(N,q)>0$.}
 \medskip

In the last section we study the singular separable solutions of $(\ref{I-0})$ written under the form
 \bel{I-9}
u(x)=u(r,\gs)=r^{-\gg_{p,q}}\gw(\gs)\qquad (r,\gs)\in \BBR_+^*\ti S^{N-1}.
\ee
Then $\gw$ satisfies the following nonlinear equation on $S^{N-1}$
 \bel{I-10}
-\Gd'\gw+\gg_{p,q}\left(N-\frac{2p+q}{p+q-1}\right)\gw-\abs\gw^{p-1}\gw\left(\gg^2_{p,q}\gw^2+\abs{\nabla'\gw}^2\right)^{\frac q2}=0,
\ee
where $\nabla'$ and $\Gd'$ are respectively the covariant gradient and the Laplace Beltrami operator on $S^{N-1}$. It is clear by integration that condition $(\ref{I-5})$ is a necessary and sufficient condition for the existence of a solution, and the constant function $\Gl_{N,p,q}$ is such a solution. We introduce a more general equation on $S^{N-1}$
 \bel{I-11}
-\Gd'\gw+\mu\gw-\abs\gw^{p-1}\gw\left(\gg^2\gw^2+\abs{\nabla'\gw}^2\right)^{\frac q2}=0,
\ee
where $\mu$ and $\gg$ are positive real numbers. We  have the following universal estimate. \medskip

 \nind{\bf Theorem E} {\it Assume $ 0\leq q<2$, $p\in\BBR$ such that $p+q-1>0$. If
   \bel{I-15}
p(N-3)+q(N-2)<N-1,
\ee
  holds,  there exists
 $c_4=c_4(N,p,q)>0$ and $a=a(N,p,q)>0$  such that for $\gg,\mu>0$, any solution $\gw$ of $(\ref{I-11})$ on $S^{N-1}$ satisfies
  \bel{I-16}
\norm\gw_{L^\infty}\leq c_4\mu^{a}\gg^{-\frac{q}{p+q-1}}.
\ee
 }

 We also give a rigidity result which shows that the solutions which are not too far from being constant are indeed constant.\medskip

 \nind{\bf Theorem F }{\it Assume $p\geq 0$ and $p+q-1>0$. Let $\gw$ be a solution of  $(\ref{I-11})$ satisfying
  \bel{I-12}
c_2^2\leq \gg^2\gw^2+\abs{\nabla'\gw}^2\leq c_1^2 ,
\ee
for some $c_1>c_2>0$ and set
  \bel{I-13}
c_*=\left\{\BA {lll} c^{p+q-1}_1\qquad\text{ if }p\geq 1\\[2mm]
c_2^{p-1}c_1^{q}\qquad\text{ if }0\leq p < 1.
\EA\right.
\ee
If
  \bel{I-14}
c_*\leq \myfrac{2(N-1+\mu)\gg^p}{q\sqrt{N-1}+2(p+q)\gg},
\ee
then $\gw$ is constant.}\medskip

In the Appendix we present many technical algebraic computations which lead to the delimitation of the regions of the $(p,q)$-plane
in which Theorem C holds. Many computations can be easily verified by using Maple. Throughout the paper $c$ denotes a generic constant depending on some parameters, specified in some cases, the value of which may change from one occurence to another.\medskip
\newpage

 The following is a representation of the different separatrix curves for $N=6$:
\begin{center}
\includegraphics[keepaspectratio, width=15cm]{curves-v5.pdf}
\end{center}

\nind{\bf Acknowledgements} This article has been prepared with the support of the collaboration
programs ECOS C14E08 and FONDECYT grant 1160540 for the three authors.

%%%%%%%%%%%%%%%%%%%%%%%%%%%%%%%%%%%%%%%%%%%%%%%%%%%%%%%%%%%%%%%%%%%%%%%%%%%%%%%%%%%%%%%%%%%%%%%%%%%%%%%%%%%%%%%%%%%%%%%%%%%%%%%%%%%%%%%%%%%%%%%%%%%%%%%%%%%%%%%%%%%%%%%%%%%%%%%%%%%%%%%%%%%%%%%%%%%%%%%%%%%%%%%%%%%%%%%%%%%%%%%%%LOCAL%ESTIMATES%%%%%%%%%%%%%%%%%%%%%%%%%%%%%%%%%%%%%%%%%%%%%%%%%%%%%%%%%%%%%%%%%%%%%%%%%%%%%%%%%%%%%%%%%%%%%%%%%%%%%%%%%%%%%%%%%%%%%%%%%%%%%%%%%%%%%%%%%%%%%%%%%%%%%%%%%%%%%%%%%%%%%%%%%%%%%%%%%%%%%%%%%%%%%%%%%%%%%%%%%%%%%%%%%%%%%%%%%%%%%%%%%%%%%%%%%%%%%%%%%%%%%%%%%%%%%%%%%%%%%%%%%%%%%%%%%%%%%%%%%%%%%%%%%%%%%%%%%%%%

\newpage

\mysection{Local estimates}
%{\color{red} Comment: I understand $N\geq 3$ in Th. 2.1. Why is it necessary in Th 2.2 and 2.3? On the other hand Th. A is for $N\geq 2$. }
\subsection{The subcritical case: Proof of Theorem A}
We show how the use of Serrin's result concerning Harnack inequality yields a blow-up estimate of any positive solution $u$ of $(\ref{I-0})$ in
a punctured domain. Assume $\overline B_1\subset\Gw$. By Brezis-Lions's result \cite{BL} there holds
  \bel{2-x1}
  u\in M^{\frac{N}{N-2}}(B_1)\,,\; \nabla u\in M^{\frac{N}{N-1}}(B_1)\,,\; u^p\abs{\nabla u}^q\in L^1(B_1),
\ee
where  $M^{r}=L^{r,\infty}$ denotes the Marcinkiewicz space or Lorentz space of index $(r,\infty)$, and there exists $\ga\geq 0$ such that
   \bel{2-x2}
-\Gd u=u^p\abs{\nabla u}^q+\ga\gd_0\quad\text{in }\CD(B_1).
\ee
We assume first $pq\neq 0$. In order to fit with Serrin's formalism, we write $u^p\abs{\nabla u}^q=\CB(u,\nabla u)$. Hence $\CB$ satisfies the estimate
   \bel{2-x3}
\abs{\CB(u,\nabla u)}\leq \abs u^{p\gth}+\abs{\nabla u}^{q\gth'}=c\abs u+d\abs{\nabla u},
\ee
where $\gth,\gth'\geq 1 $, $\frac 1\gth+\frac1{\gth'}=1$,  $c=\abs u^{p\gth-1}$, $d=\abs {\nabla u}^{q\gth'-1}$. If $\gth>\max\{1,\frac 1p\}$ and $\gth'>\max\{1,\frac 1q\}$, then $c\in M^{\frac{N}{(N-2)(p\gth-1)}}$ and  $d\in M^{\frac{N}{(N-1)(q\gth'-1)}}$. We claim that we can choose $\gth>1$ such that
 \bel{2-x5}\BA {llll}
 \myfrac{N}{(N-2)(p\gth-1)}>\myfrac{N}{2}\quad\text{and }\;\; \myfrac{N}{(N-1)(q\gth'-1)}>N.
\EA \ee
These inequalities are respectively equivalent to
 \bel{2-x5'}\BA {llll}
\gth<\myfrac{N}{p(N-2)}\quad\text{and }\;\; \gth'<\myfrac{N}{q(N-1)},
\EA \ee
which is clearly possible from $(\ref{I-6})$ by taking $\gth=\frac{N(1-\varepsilon )}{p(N-2)}$ for $\varepsilon >0$ small enough.
%%%%%%\Longleftrightarrow \\[4mm]%%%%%%%%%%%%%%%%%%%%%%%%%%%%%%%%%%%%%%%%%%%%%%%\Longleftrightarrow %%%%%%%%%%%%%%%%%%%%%%%%%%%%%%%%%%%%%%%%%
Because $M^r(B_1)\subset_{\!\!>}L^{r-\gd}(B_1)$ for any $\gd>0$, we infer that $c\in L^{\frac{N}{2}+\gd}(B_1)$ and $c\in L^{N+\gd}(B_1)$ and $u$ verifies Harnack inequality in $B_{1}\setminus\{0\}$ by \cite[Th 5]{Se1}. This implies
 \bel{2-x6}\BA {llll}\displaystyle
\max_{\abs x=r}u(x)\leq K\min_{\abs x=r}u(x)\quad\forall r\in (0,\tfrac12]\quad\mbox{for some $K>0$}.
\EA \ee
The spherical average $\bar u$ of $u$ on $\{x:\abs x=r\}$ is superharmonic. Hence there exists some $m\geq 0$ such that
 \bel{2-x7}\BA {llll}\displaystyle
\bar u(r)\leq m r^{2-N}.
\EA \ee
Combined with $(\ref{2-x6})$ it yields $u(x)\leq Km\abs x^{2-N}$. The estimate on the gradient is standard, see eg \cite[Lemma 3.3.2]{Vebook2}.\qeda\medskip

\nind\Remark Estimate $(\ref{I-7})$ is not universal since the constant $K$ in $(\ref{2-x6})$ depends on the norms of $c$ and $d$ which could depend not only on $N$, $p$, $q$, but on the solution itself.\medskip

The following result is not new, except in the case $p=0$, $q=1$. It was proved in \cite[Th 15.1]{MiPo} for $p+q-1>0$ and extended to quasilinear operators by simulating  the change of unknown $u=v^b$. It was used in order to derive a priori estimates \cite{BVP}. Later on it was extended to more general operators in \cite{Fil} where the new cases $p+q-1<0$ and $p+q-1=0$ with $p>0$ were considered using a delicate proof. We give here a very simple  but general proof of all these results. Furthermore our method highlights the role of the change of unknown function, which foreshadows the method used in Theorem B. It is also extendable to very general quasilinear operators such as
\bel{2-x10-}
-div A(x,u,\nabla u)\geq B(x,u,\nabla u)\quad\text{in }\;\BBR^N,
\ee
under the assumptions that $\langle A(x,r,\xi),\xi\rangle \geq \abs\xi^m$, $B(x,r,\xi)r\leq c\abs r^{p}\abs\xi^q$, and under the corresponding subcritical condition $p(N-m)+q(N-1)<N(m-1)$.

\bth{nonexis} Assume $N\geq 2$, $p$ and $q$ are nonnegative and $(\ref{I-6})$ holds. Then the only positive functions $u\in C^1(\BBR^N)$ satisfying
 \bel{2-x8'}\BA {llll}\displaystyle
-\Gd u\geq u^p\abs{\nabla u}^q,
\EA \ee
in $\BBR^N$ are the constants.
 \es
\Proof Assume $u$ is such a solution. For $p+q\neq 1$, we set $u=v^{b}$ with $b(b-1)>0$.
$$-b\Gd v\geq b(b-1)\myfrac{\abs{\nabla v}^2}{v}+|b|^{q}v^s\abs{\nabla v}^q,
$$
with
$$s=1-q+b(p+q-1).
$$
If $s>0$, then from H\"older's inequality,
$$\abs{\nabla v}^{\frac{2s+q}{s+1}}=\left(\myfrac{\abs{\nabla v}^2}{v}\right)^{\frac{s}{s+1}}v^{\frac{s}{s+1}}\abs{\nabla v}^{\frac{q}{s+1}}
\leq\gd^{\frac{s+1}{s}}\myfrac{\abs{\nabla v}^2}{v}+ \gd^{-1-s}v^s\abs{\nabla v}^{q},$$
for any $\gd>0$. Hence, by chosing $\gd$, we see that there exists $c>0$ such that
\bel{2-x8}-b\Gd v\geq c\abs{\nabla v}^\gw\quad\text{where }\;\gw=\myfrac{2s+q}{s+1}=\myfrac{2-q+2b(p+q-1)}{2-q+b(p+q-1)}.
\ee

(i) In the case $p+q-1>0$, we take $b=1+\varepsilon $, for $\varepsilon >0$. Then $s=p+\varepsilon (p+q-1)>0$, and $s>1-q$, thus $\gw>1$. From assumption $(\ref{I-6})$ we can take $\varepsilon >0$ small enough such that $(N-2)s+(N-1)q<N$, which is equivalent to $\gw<\frac{N}{N-1}$.\smallskip

(ii) If $p+q-1<0$, hence $0\leq q<1$, we take $b=-\varepsilon $, for $\varepsilon >0$. Then

\bel{2-x9}-|b|\Gd v+ c\abs{\nabla v}^\gw\leq 0,\quad\text{where }\; c>0.
\ee
and $s=1-q-\varepsilon  (p+q-1))>1-q>0$. hence $\gw=1-\frac{p+q-1}{2-q}>1$ and we can choose $\varepsilon $ small enough such that $\gw<\frac{N}{N-1}$.\smallskip

(iii) If $p+q-1=0$, we set $u=e^v$ (where $v$ is a signed function) and derive that for any $\tilde \gw\in (q,2)$ one can find $\tilde c>0$ such that
$$-\Gd v\geq \abs{\nabla v}^2+e^{p+q-1}\abs{\nabla v}^q\geq \tilde c\abs{\nabla v}^{\tilde \gw}.
$$
We can take in particular $\tilde \gw<\frac{N}{N-1}$.\smallskip

Finally, let $R>0$ and
$\gz\in C_0^\infty(\BBR^N)$, with values in $[0,1]$, such that $\gz=1$ on $B_{\frac R2}$, $\gz=0$ on $B_R^c$ and $\abs{\nabla\gz}\leq 2R^{-1}$.
Then, in each of the three cases,  there exists $c_1>0$ such that
$$\BA {lll}\myint{B_R}{}\gz^{\gw'}\abs{\nabla v}^\gw dx\leq c_1\left|\myint{B_R}{}\langle\nabla v,\nabla\gz\rangle\gz^{\gw'-1}dx\right|
\\[4mm]
\phantom{\myint{B_R}{}\gz^{\gw'}\abs{\nabla v}^\gw dx}
\leq \myfrac{1}{2}\myint{B_R}{}\gz^{\gw'}\abs{\nabla v}^\gw dx+c_2\myint{B_R}{}\abs{\nabla \gz}^{\gw'}dx,
\EA$$
where $\gw'$ is the H\"older conjugate of $\gw$, from which follows
$$\BA {lll}\myint{B_{\frac R2}}{}\abs{\nabla v}^\gw dx\leq 2C\myint{B_R}{}\abs{\nabla \gz}^{\gw'}dx\leq C'R^{N-\gw'}.
\EA$$
which implies the claim since $\gw'>N$.\qeda \medskip

\nind\Remark If the inequality $(\ref{2-x8})$ is considered in an exterior domain, the situation differs according to whether $0\leq q<1$ or $1\leq q<2$. The case $q=0$ is known for a long time. The following result is proved in \cite{BGMQ}: {\it
 Assume $0<q<1$, $f\in C(\BBR_+)$ is positive on $(\BBR_+)^*$ and satisfies
\bel{2-x10}\myint{0}{1}\myfrac{f(s)}{s^{\gth}}ds=\infty\quad\text {with }\;\gth=\myfrac{(2-q)(N-1)}{N-2}.
\ee
Then every $C^1$ positive function $u$ satisfying
\bel{2-x11}-\Gd u\geq f(u)|\nabla u|^q\quad\text {in }\;\Gw=\BBR^N\setminus K
\ee
for some compact set $K$ becomes eventually constant.}
%%%%%%%%%%%%%%%%%%%%%%%%%%%%%%%%%%%%%%%%%%%%%%%%%%%%%%%%%%%%%%%%%%%%%%%%%%%%%%%%%%%%EXTERIOR%%DOMAIN%%%%%%%%%%%%%%%%%%%%%%%%%%%%%%%%%%%%%%%%%%%%%%%%%%%%%%%%%%%%%%%%%%%%%%%%%%%%%%%%%%%%%%%%%%%%%%%%%%%%%%%%%%%%%%%%%%%%%%%%%%%%%%%%%%%%%%%%%%%%%%%%%%%%%%%%%%%%%%%%

%%%%%%%%%%%%%%%%%%%%%%%%%%%%%%%%%%%%%%%%%%%%%%%%%%%%%%%%%%%%%%%%%%%%%%%%%%%%%%%%%%%%%%%%%%%%%%%%%%%%%%%%%%%%%%%%%%%%%%%%%%%%%%%%%%%%%%%%%%%%%%%%%%%%%%%%%%%%%%%%%%%%%%%%%%%%%%%%%%%%%%%%%%%%%%%%%%%%%%%%
%%%%%%%%%%%%%%%%%%%%%%BERNSTEIN%%The general case%%%%%%%%%%%%%%%%%%%%%%%%%%%%%%%%%%%%%%%%%%%%%%%%%%%%%%%%%%%%%%%%%%%%%%%%%%%%%%%%%%%%%%%%%%%%%%%%%%%%%%%%%%%%%%%%%%%%%%%%%%%%%%%%%%%%%%%%%%%%%%%%%%%%%%%%%%%%%%%%%%%%%%%%

\subsection{Proof of Theorem B}

The next result will be useful in the sequel.
\blemma{Oss} Let $q>1$ and $R>0$. Assume $\gu$ is continuous and nonnegative on $\overline B_R$ and $C^1$ on the set $\CU_+=\{x\in B_R:\gu(x)>0\}$. If  $\gu$ satisfies for some real number $a$,
  \bel{2-00}
-\Gd\gu+\gu^q\leq a\myfrac{\abs{\nabla \gu}^2}{\gu}
\ee
 on each connected component of $\CU_+$, there holds
  \bel{2-01}
\gu(0)\leq c_{N,q,a}R^{-\frac{2}{q-1}}.
\ee
\es
\Proof We can always suppose $a>0$ and set $\gu=W^\ga$ for some $\ga>0$  to be defined. Then
$$-\Gd W+(1-\ga)\myfrac{\abs{\nabla W}^2}{W}+ \myfrac{1}{\ga}W^{\ga(q-1)+1}\leq a\ga^2\myfrac{\abs{\nabla W}^2}{W}.
$$
If we choose $1-\ga\geq a\ga^2$, or equivalently $0<\ga\leq \frac{1}{a+1}$ we derive
$$-\Gd W+\myfrac{1}{\ga}W^{\ga(q-1)+1}\leq 0.
$$
on each connected component of $\CU_+$. A standard computation shows that there exists $c_{N,\ga,q}>0$ such that the function
$$x\mapsto \psi(x):=\myfrac{c_{N,\ga,q}(R^2\ga)^{\frac{1}{\ga (q-1)}}}{\left(R^2-|x|^2\right)^{\frac{2}{\ga (q-1)}}}
$$
satisfies
$$-\Gd \psi+\myfrac{1}{\ga}\psi^{\ga(q-1)+1}\geq 0 \qquad\text{in }\; B_R.
$$
Assume that there exists a connected component $G$ of $\{x\in B_R:W(x)>\psi(x) \}$, then $G\subset\CU_+$. The function
$\phi=W-\psi$ is subharmonic and continuous in $\overline G$ and vanishes on $\prt G$. Hence $\phi\leq 0$, contradiction. Hence
$G=\emptyset$ and $W\leq \psi $ in $\CU_+$. Since $W\equiv 0$ in $B_R\setminus\CU_+$ we derive that $W\leq \psi $ in $B_R$. Therefore
  \bel{2-02}
W(0)\leq \myfrac{c_{N,\ga,q}\ga^{\frac{1}{\ga (q-1)}}}{R^{\frac{2}{\ga (q-1)}}}
\ee
which leads to $(\ref{2-01})$.\qeda\medskip

%%%%%%%%%%%%%%%%%%%%%%%%%%%%%%%%%%%%%%%%%%%%%%%%%%%%%%%%%%%%%%%%%%%%%%%%%%%%%%%%%%%%%%%%%%%%%%%%%%%%%%%%%%%%%%%%%%%%%%%%%%%%%%%%%%%%%%%%%%%%%%%%%%%%%%%%%%%%%PROOF OF Theorem B%%%%%%%%%%%%%%%%%%%%%%%%%%%%%%%%%%%%%%%%%%%%%%%%%%%%%%%%%%%%%%%%%%%%%%%%%%%%%%%%%%%%%%%%%%%%%%%%%%%%%%%%%%%%%%%%%%%%%%%%%%%%%%%%%%%%%%%%%%%%%%%%%%%%%%%%%%%%%%%%%%%%%%%%%%%%%%%%%%%

\nind{\it Proof of Theorem B}: In any open subset $\CU$ of $B_R$ where $\abs{\nabla u}>0$ the function $u$ is $C^\infty$ and the next computations are justified.\smallskip

\nind{\it Step 1: Transformation of the equation}. Set $u=v^{-\gb}$ where $\gb$ is a nonzero real number to be chosen. Then
\bel{2-1}\BA {lll}
\Gd v=(1+\gb)\myfrac{\abs {\nabla v}^2}{v}+\abs\gb^{q-2}\gb v^{1-q-\gb(p+q-1)}\abs {\nabla v}^q\\[4mm]
\phantom{\Gd v}
=(1+\gb)\myfrac{z}{v}+\abs\gb^{q-2}\gb v^sz^{\frac q2},
\EA\ee
if we denote $z=\abs {\nabla v}^2$ and $s=1-q-\gb(p+q-1)$. We recall that
$$\frac 12\Gd\abs {\nabla v}^2=\abs{D^2 v}^2+\langle\nabla\Gd v,\nabla v\rangle.$$
Since
$$\abs{D^2 v}^2\geq \frac 1N\left(\Gd v\right)^2= \frac 1N\left((1+\gb)^2\myfrac{z^2}{v^2}+\gb^{2(q-1)} v^{2s}z^{q}+
2(1+\gb)\abs\gb^{q-2}\gb v^{s-1}z^{\frac q2+1}\right),
$$
we get
 \bel{2-2}\BA {lll}
 \myfrac 12\Gd z\geq \myfrac 1N\left((1+\gb)^2\myfrac{z^2}{v^2}+\gb^{2(q-1)} v^{2s}z^{q}+
2(1+\gb)\abs\gb^{q-2}\gb v^{s-1}z^{\frac q2+1}\right)-(1+\gb)\myfrac{z^2}{v^2}\\[4mm]
\phantom{ \myfrac 12\Gd z}
+s\abs\gb^{q-2}\gb v^{s-1}z^{\frac q2+1}+(1+\gb)\myfrac{\langle\nabla z,\nabla v\rangle}{v}
+\myfrac{q}{2}\abs\gb^{q-2}\gb  v^sz^{\frac q2-1}\langle\nabla z,\nabla v\rangle,
 \EA\ee
 which can be re-written as
 \bel{2-3}\BA {lll}
 - \myfrac 12\Gd z+\left(\myfrac {(1+\gb)^2}N-(1+\gb)\right)\myfrac{z^2}{v^2}+\myfrac 1N\gb^{2(q-1)} v^{2s}z^{q}
 +\left(\myfrac {2(1+\gb)}{N}+s\right)\abs\gb^{q-2}\gb v^{s-1}z^{\frac q2+1}\\[4mm]
 \phantom{---------------}
 +(1+\gb)\myfrac{\langle\nabla z,\nabla v\rangle}{v}+\myfrac{q}{2}\abs\gb^{q-2}\gb  v^sz^{\frac q2-1}\langle\nabla z,\nabla v\rangle\leq 0.
 \EA\ee
Next we set $z=v^{-\gl}Y $ for some parameter $\gl$, then
$$-\Gd z=\gl v^{-\gl-1}Y\Gd v-\gl(\gl+1)v^{-2\gl-2}Y^2+2\gl v^{-\gl-1}\langle\nabla Y,\nabla v\rangle-v^{-\gl}\Gd Y.
$$
Since $\Gd v=(1+\gb)\frac zv+\abs\gb^{q-2}\gb v^{s}z^{\frac q2}=(1+\gb)v^{-\gl-1}Y+\abs\gb^{q-2}\gb v^{s-\frac{\gl q}{2}}Y^{\frac q2}$ we use
$$\frac{z^2}{v^2}=v^{-2\gl-2}Y^2\,,\;v^{2s}z^q=v^{2s-\gl q}Y^q\,\text{ and }\;v^{s-1}z^{1+\frac q2}=v^{s-1-\gl-\frac{\gl q}{2}}Y^{\frac q2+1},$$
and get
$$-\Gd z=\gl(\gb-\gl)v^{-2-2\gl}Y^2+\gl\abs\gb^{q-2}\gb v^{s-1-\gl-\frac{\gl q}{2}}Y^{\frac q2+1}+2\gl v^{-\gl-1}\langle\nabla Y,\nabla v\rangle
-v^{-\gl}\Gd Y.
$$
 Reporting into $(\ref{2-3})$ yields
$$\BA {lll}\displaystyle
0\geq\myfrac{\gl}{2}(\gb-\gl)v^{-2-2\gl}Y^2+\myfrac{\gl}{2}\abs\gb^{q-2}\gb v^{s-1-\gl-\frac{\gl q}{2}}Y^{\frac q2+1}+\gl v^{-\gl-1}\langle\nabla Y,\nabla v\rangle-\myfrac{1}{2}v^{-\gl}\Gd Y\\[4mm]\phantom{0\geq}
+\left(\myfrac {(1+\gb)^2}N-(1+\gb)\right)v^{-2\gl-2}Y^2
+\left(\myfrac {2(1+\gb)}{N}+s\right)\abs\gb^{q-2}\gb v^{s-1-\gl-\frac{\gl q}{2}}Y^{\frac q2+1}\\[4mm]\phantom{0\geq}
+\myfrac 1N\gb^{2(q-1)} v^{2s-\gl q}Y^{q}+(1+\gb)\left(v^{-\gl-1}\langle\nabla Y,\nabla v\rangle-\gl v^{-2\gl-2}Y^2\right)\\[4mm]\phantom{0\geq}
+\myfrac q2\abs\gb^{q-2}\gb\left(v^{s-\frac{\gl q}{2}}Y^{\frac{q}{2}-1}\langle\nabla Y,\nabla v\rangle-\gl v^{s-1-\gl-\frac{\gl q}{2}}Y^{1+\frac q2}\right),
\EA$$
which can be re-written under the form
$$\BA {lll}\displaystyle
-\myfrac{1}{2}v^{-\gl}\Gd Y+\left(\myfrac {(1+\gb)^2}N-(1+\gb)-\myfrac{\gl}{2}(\gl+\gb+2)\right) v^{-2\gl-2}Y^2\\[4mm]\phantom{\leq--}
+\left(\myfrac{2(1+\gb)}{N}+s-\myfrac{\gl(q-1)}{2}\right)\abs\gb^{q-2}\gb v^{s-1-\gl-\frac{\gl q}{2}}Y^{\frac q2+1}+\myfrac 1N\gb^{2(q-1)} v^{2s-\gl q}Y^{q}\\[4mm]\phantom{\leq--}
\leq -\left(\myfrac{q}{2}\abs\gb^{q-2}\gb v^{s-\frac{\gl q}{2}}Y^{\frac q2-1}+(\gl+1)v^{-\gl-1}\right)\langle\nabla Y,\nabla v\rangle.
\EA$$
Multiplying this relation by $v^\gl$ yields
$$\BA {lll}\displaystyle
-\myfrac{1}{2}\Gd Y+\left(\myfrac {(1+\gb)^2}N-(1+\gb)-\myfrac{\gl}{2}(\gl+\gb+2)\right) v^{-\gl-2}Y^2\\[4mm]\phantom{\leq--}
+\left(\myfrac{2(1+\gb)}{N}+s-\myfrac{\gl(q-1)}{2}\right)\abs\gb^{q-2}\gb v^{s-1-\frac{\gl q}{2}}Y^{\frac q2+1}+\myfrac 1N\gb^{2(q-1)} v^{2s-\gl q+\gl}Y^{q}\\[4mm]\phantom{\leq--}
\leq -\left(\myfrac{q}{2}\abs\gb^{q-2}\gb v^{s-\frac{\gl q}{2}+\gl}Y^{\frac q2-1}+\myfrac{\gl+1}{v}\right)\langle\nabla Y,\nabla v\rangle.
\EA$$\smallskip

\nind {\it Step 2: Estimate on $Y$}.
Let $\varepsilon _0\in (0,1)$. For any $\varepsilon >0$ one has
$$\abs{\myfrac{\langle\nabla Y,\nabla v\rangle}{v}}
\leq \myfrac{1}{4\varepsilon }\myfrac{\abs{\nabla Y}^2}{Y}+\varepsilon  v^{-\gl-2} Y^2.
$$
Taking $\varepsilon =\frac{\varepsilon _0}{\abs{\gl+1}}$, we get
$$\abs{(\gl+1)\myfrac{\langle\nabla Y,\nabla v\rangle}{v}}\leq \myfrac{(\gl+1)^2}{4\varepsilon _0}\myfrac{\abs{\nabla Y}^2}{Y}+\varepsilon _0 v^{-\gl-2} Y^2.
$$
In the same way, with $\varepsilon' =\frac{2\varepsilon _0}{q\abs\gb^{q-1}}$,
$$\abs{v^{s-\frac{\gl q}{2}+\gl}Y^{\frac q2-1}\langle\nabla Y,\nabla v\rangle}=
v^{s-\frac{\gl q}{2}+\gl}Y^{\frac q2-1}Y^{\frac{q-1}{2}}\myfrac{\abs{\nabla Y}}{\sqrt Y}\leq \varepsilon 'v^{2s-\gl q+\gl}Y^q+\myfrac{1}{4\varepsilon '}
\myfrac{\abs{\nabla Y}^2}{Y},
$$
and
$$-\frac q2\abs\gb^{q-2}\gb v^{s-\frac{\gl q}{2}+\gl}Y^{\frac q2-1}\langle\nabla Y,\nabla v\rangle\leq\frac{\varepsilon 'q\abs\gb^{q-1}}{2}
v^{2s-\gl q+\gl}Y^q+\frac{q\abs\gb^{q-1}}{8\varepsilon '}\frac{\abs{\nabla Y}^2}{Y}.
$$
We infer
 \bel{2-4}\BA {lll}\displaystyle
 -\frac 12\Gd Y+\left(\frac{(1+\gb)^2}{N}-(1+\gb)-\frac{\gl(\gb+\gl+2-\varepsilon _0)}{2}\right)v^{-\gl-2}Y^2\\[4mm]\phantom{--------}\displaystyle
 +\left(\frac{2(1+\gb)}{N}+s-\frac{\gl(q-1)}{2}\right)\abs\gb^{q-2}\gb v^{s-1-\frac{\gl q}{2}}Y^{1+\frac{q}{2}}\\[4mm]\phantom{--------}\displaystyle+\left(\frac{\gb^{2(q-1)}}{N}-\varepsilon _0\right)v^{2s-\gl q+\gl}Y^q\leq C(\varepsilon _0)\frac{\abs{\nabla Y}^2}{Y},
\EA\ee
with $C(\varepsilon _0)=\left(\frac{(\gl+1)^2}{4}+\frac{q\gb^{2(q-1)}}{16}\right)\frac{1}{\varepsilon _0}$. Next we put
%%%%%%%%%%%%%%%%%%%%%%%%%%%%%%%%%%%%%%%%%%%%%%%%%%%%%
$$\BA {lll}\displaystyle H=\left(\frac{(1+\gb)^2}{N}-(1+\gb)-\frac{\gl(\gb+\gl+2-\varepsilon _0)}{2}\right)v^{-\gl-2}Y^2\\[4mm]\phantom{}\displaystyle
+\left(\frac{2(1+\gb)}{N}+s-\frac{\gl(q-1)}{2}\right)\abs\gb^{q-2}\gb v^{s-1-\frac{\gl q}{2}}Y^{1+\frac{q}{2}}
+\left(\frac{\gb^{2(q-1)}}{N}-\varepsilon _0\right)v^{2s-\gl q+\gl}Y^q,
\EA$$
and consider the trinomial
$$\BA {lll}\displaystyle
{\bf T}_{\varepsilon _0}(t)=\left(\frac{\gb^{2(q-1)}}{N}-\varepsilon _0\right)t^2+\left(\frac{2(1+\gb)}{N}+s-\frac{\gl(q-1)}{2}\right)\abs\gb^{q-2}\gb t\\[4mm]
\phantom{{\bf T}_{\varepsilon _0}(t)}\displaystyle
+\left(\frac{(1+\gb)^2}{N}-(1+\gb)-\frac{\gl(\gb+\gl+2-\varepsilon _0)}{2}\right).
\EA$$
If its discriminant is negative there exists $\ga=\ga(N,p,q,\gb,\gl,\varepsilon _0)>0$ such that ${\bf T}_{\varepsilon _0}(t)\geq \ga(t^2+1)$, hence
 \bel{2-5}\BA {lll}\displaystyle
H\geq\ga\left(v^{-\gl-2}Y^2+v^{2s-\gl q+\gl}Y^q\right).
\EA\ee
Assuming $\gl\neq -2$, we introduce
 \bel{2-S}S=\myfrac{2s-\gl q+\gl}{\gl+2}=1-q-\myfrac{2\gb(p+q-1)}{\gl+2},
\ee
then, if $S>\max\{0,1-q\}$, we have $\frac{2S+q}{S+1}>1$ and
$$Y^{\frac{2S+q}{S+1}}=\left(\myfrac{Y^2}{v^{\gl+2}}\right)^{\frac{S}{S+1}}v^{\frac{(\gl+2)S}{S+1}}Y^{\frac{q}{S+1}}\leq \myfrac{Y^2}{v^{\gl+2}}+
v^{(\gl+2)S}Y^q= \myfrac{Y^2}{v^{\gl+2}}+v^{2s-\gl q+\gl}Y^q.
$$
From this we infer the key inequality
 \bel{2-6}\BA {lll}\displaystyle
-\Gd Y+2\ga Y^{\frac{2S+q}{S+1}}\leq 2C(\varepsilon _0)\frac{\abs{\nabla Y}^2}{Y}.
\EA\ee
Using \rlemma{Oss}, we derive
 $$
Y(0)\leq cR^{-\frac{2(S+1)}{S+q-1}}=cR^{-\frac{2(s+1)-\gl(q-1)}{s+q-1}}=cR^{-2+\frac{(2+\gl)(2-q)}{\gb(p+q-1)}},
 $$
 from which follows
  \bel{2-7}
  \abs{\nabla u^{-\frac{2+\gl}{2\gb}}(0)}\leq \myfrac{\abs{2+\gl}\sqrt c}{2}R^{-1+\frac{(2+\gl)(2-q)}{2\gb(p+q-1)}}.
\ee
 Therefore, $(\ref{I-8})$ in Theorem B will follow with  $a=-\frac{\gl+2}{2\gb}$ and $a$ will be positive  from $(\ref{2-S})$ if we can choose $\gl\not=-2$ and $\gb\not=0$ so that $S>\max\{0,1-q\}$. In what follows we shall see that under the assumptions of Theorem B we can always choose such $\gb$, $\gl$.
\smallskip

\nind {\it Step 3: Study of the trinomial ${\bf T}_{\varepsilon _0}$}. The discriminant of the trinomial ${\bf T}_{\varepsilon _0}$ is a polynomial in its coefficients. Hence it is sufficient to prove that the discriminant of ${\bf T}_0$ is negative to derive that the same property holds for ${\bf T}_{\varepsilon _0}$ for $\varepsilon _0$ small enough.
If
$$\BA {lll}\displaystyle
{\bf T}_0(t)=\frac{\gb^{2(q-1)}}{N}t^2+\left(\frac{2(1+\gb)}{N}+s-\frac{\gl(q-1)}{2}\right)\abs\gb^{q-2}\gb t\\[4mm]
\phantom{{\bf T}_{\varepsilon _0}(t)_0--------}\displaystyle
+\left(\frac{(1+\gb)^2}{N}-(1+\gb)-\frac{\gl(\gb+\gl+2)}{2}\right),
\EA$$
its discriminant $D$ verifies
$$\BA {lll}\displaystyle
\gb^{2(1-q)}D=\left(\frac{2(1+\gb)}{N}+s-\frac{\gl(q-1)}{2}\right)^2-\frac{4}{N}\left(\frac{(1+\gb)^2}{N}-(1+\gb)-\frac{\gl(\gb+\gl+2)}{2}\right).
\EA$$
Using $\gb+1=\frac{p-s}{p+q-1}$ we obtain
$$\gb^{2(1-q)}D=\left(s-\gl\frac{q-1}{2}\right)^2+\frac{4(p-s)}{N(p+q-1)}\left(s+1+\gl\frac{2-q}{2}\right)+\frac{2\gl(\gl+1)}{N}.
$$
Since $S=\frac{2s+\gl(1-q)}{\gl+2}$, $s-\gl\frac{q-1}{2}=\frac{(\gl+2)S}2$, hence
$$\BA{lll}\displaystyle
\gb^{2(1-q)}D=\frac{(\gl+2)^2S^2}{4}+\frac{4}{N(p+q-1)}\left[-\frac{(\gl+2)^2S^2}{4}+\left(p-1-\frac{\gl q}{2}\right)\frac{(\gl+2)S}{2}\right.
\\[4mm]\phantom{--------------}\displaystyle
\left.+\left(p-\frac{\gl(q-1)}{2}\right)\left(\frac{\gl+2}{2}\right)\right]+\frac{2\gl(\gl+1)}{N}.
\EA$$

Set $Q=p+q-1$ and $D_1=NQ\gb^{2(1-q)}D$, then
$$D_1=(\gl+2)^2\left(\myfrac{NQ}{4}-1\right)S^2+2\left(p-1-\myfrac{\gl q}{2}\right)(\gl+2)S+\tilde L,
$$
where miraculously,
$$\BA {lll}\tilde L=(2p+\gl(1-q))(\gl+2)+2\gl(\gl+1)Q\\\phantom{\tilde L}
=Q\gl^2+p(\gl+2)^2>0.
\EA$$
So we require that $\gl+2\neq 0$ and set $\ell=\frac{\gl}{\gl+2}$; hence $\ell\neq 1$ and $\gl+2=\frac{2}{1-\ell}$. We obtain
  \bel{2-8}\BA {lll}\displaystyle
  D_2(S,\ell):=\myfrac{D_1}{(\gl+2)^2}=\left(\myfrac{NQ}{4}-1\right)S^2+(p-1-Q\ell)S+Q\ell^2+p,
  \EA\ee
and we look for $\ell\neq 1$ and $S>(1-q)_+$ such that $D_2(S,\ell)<0$. \\
  	
 	\nind We can write
 $$D_2(S,\ell)=Q\left(\ell-\frac{S}{2}\right)^2+ \left(\myfrac{(N-1)Q}{4}-1\right)S^2-(1-p)S+p:=Q\left(\ell-\frac{S}{2}\right)^2+\CT(S)$$
 and we need to find $\ell\neq 1$ and $S>(1-q)^+$ such that $D_2(S,\ell)<0$.

  \nind (i) We first assume $Q<\frac{4}{N-1}$. We fix $\ell=\frac{S}{2}$. As the coefficient of $S^2$ in
  $\CT(S)$ is negative, we can choose $S$ large enough and $S>2>(1-q)_+$ so  that $\CT(S)<0$ and $\ell>1$.

 \nind (ii) Let $0\leq p< 1$. and note that our assumption $Q<\frac{(p+1)^2}{(N-1)p}$ is equivalent to
 \bel{2-10m}\BA {lll}\displaystyle
 d:=(p-1)^2-p\left((N-1)Q-4\right)>0
 \EA\ee
that is, the discriminant of $\CT(S)$ is positive. We are only concerned with the case $Q\geq \frac{4}{N-1}$.
Assume first that $Q= \frac{4}{N-1}$.  Then it suffices to choose $\ell=\frac{S}{2}$ with $S>\max\{\frac{p}{1-p},2\}$ in order to achieve our goal.

Finally assume $Q> \frac{4}{N-1}$. Then, as $p<1$, the trinomial $\CT(S)$ has two positive roots $S_1$ and $S_2$. We will see next that $$S=\frac{S_1+S_2}{2}=\myfrac{2(1-p)}{(N-1)Q-4}>0$$
and $\ell=\frac{S}{2}+\varepsilon $ where $\varepsilon =0$ if $S\neq2$,  and
$$0<\varepsilon ^2<\frac{d}{Q((N-1)Q-4)}\quad\mbox{if}\quad S=2$$ satisfy all the requirements. Indeed, from $(\ref{2-10m})$ we have
$$S+q-1=\myfrac{2(1-p)}{(N-1)Q-4}+q-1>\myfrac{2p}{1-p}+q-1=\myfrac{p(2-q)+Q}{1-p}>0,$$ hence $S>(1-q)^+$,  and since $d$ is the discriminant of the trinomial $\CT(S)$ we have that
$$D_2(S,\ell)=Q\varepsilon ^2-\frac{d}{(N-1)Q-4}<0.$$

At end we choose $\ell\neq 1$ close enough to $\frac{S_0}{2}$ so that all the conditions are satisfied (always with $a=-\frac{\gl+2}{2\gb}>0$). \qeda
\medskip

\nind\Remark In the case $Q<\frac 4{N-1}$, we have fixed some $\ell>1$ so that $\gl+2<0$ and then $\gb>0$. If $Q<\frac 4{N}$, a much simpler possible choice is $\gl=\ell=0$ so that $\gb<0$.

%%%%%%%%%TEXTCOLOR-RED-%%%%%%%%%%%%%%%%%%%%%%%%%%%%%%%%%%%%%%%%%%%%%%%%%%%%%%%%%%%%%%%%%%%%%%%%%%%%%%%%%%%%%%%%%%%%%%%%%%%%%%%%%%%%%%%%%%%%%
\medskip

%%%%%%%%%%%%%%%%%%%%%%%%%%%%%%%%%%%%%%%%%%%%%%%%%%%%%%%%%%%%%%%%%%%%%%%%%%%%%%%%%%%%%%%%COROLLARY%%%%%%%%%%%%%%%%%%%%%%%%%%%%%%%%%%%%%%%%%%%%%%%%%%%%%%%%%%%%%%%%%%%%%%%%%%%%%%%%%%%%%%%%%%%%%%%%%%%%%%%%%%%%%%%%%%%%%%%%%%%%%%%%%%%%%%%%%%%%%%%%%%%%%%%%%%%%%%%%%%%%%%%%%%%%%

\nind{\it Proof of Corollary B-2}. Since $\prt\Gw$ is smooth, there exists $d_0>0$ such that for any $z\in \Gw$ verifying $\dist (z,\prt\Gw)\leq d_0$, there exists a unique $\gz_z\in\prt\Gw$ such that $\dist (z,\prt\Gw)=\abs{z-\gz_z}$. If $\dist (z_0,\prt\Gw)=d_0$, we denote by ${\bf n}_{\gz_{z_0}}$ the normal inward unit vector to $\prt\Gw$ at $\gz_{z_0}$ and we set $x_t=t{\bf n}_{\gz_{z_0}}$, $0<t\leq d_0$. By  $(\ref{I-8})$,
$$\abs{\nabla u^a(x_t)}\leq c_1t^{-1-a\frac{2-q}{p+q-1}},
$$
hence
$$\abs{u^a(x_t)-u^a(z_0)}\leq c_1\myint{t}{d_0}s^{-1-a\frac{2-q}{p+q-1}}ds.
$$
This implies
$$u^a(x_t)\leq u^a(z_0)+\frac{c_1(p+q-1)}{a(2-q)}t^{-a\frac{2-q}{p+q-1}}=u^a(z_0)+\frac{c_1(p+q-1)}{a(2-q)}\left(\dist(x_t,\prt\Gw)\right)^{-a\frac{2-q}{p+q-1}}.
$$
If $a\geq 1$ it yields
$$u(x_t)\leq u(z_0)+\frac{c_1(p+q-1)}{(2-q)}\left(\dist(x_t,\prt\Gw)\right)^{-\frac{2-q}{p+q-1}},
$$
while, if $0<a<1$ we can only obtain
$$u(x_t)\leq c_2\left((u(z_0)+\frac{c_1(p+q-1)}{(2-q)}\left(\dist(x_t,\prt\Gw)\right)^{-\frac{2-q}{p+q-1}}\right).
$$
In any case we derive $(\ref{I-8-1})$.\qeda
%%%%%%%%%%%%%%%%%%%%%%%%%%%%%%%%%%%%%%%%%%%%%%%%%%%%%%%%%%%%%%%%%%%%%%%%%%%%%%%%%%%%%%%%%%%%%%%%%%%%%%%%%%%%%%%%%%%%%%%%%%%%%%%%%%%%%%%%%%%%%%%%%%%%%%%%%%%%%GLOBAL%%%%%%%%%%%%%%%%%%%%%%%%%%%%%%%%%%%%%%%%%%%%%%%%%%%%%%%%%%%%%%%%%%%%%%%%%%%%%%%%%%%%%%%%%%%%%%%%%%%%%%%%%%%%%%%%%%%%%%%%%%%%%%%%%%%%%%%%%%%%%%%%%%%%%%%%%%%%%%%%%%%%%%%%%%%%%%
%%%%%%%%%%%%%%%%%%%%%%%%%%%%%%%%%%%%%%%%%%%%%%%%%%%%%%%%%%%%%%%%%%%%%%%%%%%%%%%%%%%%%%%%%%%%%%%%%%%%%%%%%%%%%%%%%%%%%%%%%%%%%%%%%%%%%%%%%%%%%%%%%%%%%%%%%%%%%SOLUTIONS%%%%%%%%%%%%%%%%%%%%%%%%%%%%%%%%%%%%%%%%%%%%%%%%%%%%%%%%%%%%%%%%%%%%%%%%%%%%%%%%%%%%%%%%%%%%%%%%%%%%%%%%%%%%%%%%%%%%%%%%%%%%%%%%%%%%%%%%%%%%%%%%%%%%%%%%%%%%%%%%%%%%%%%%%%%%%%%%%%%%%%%%%%%%%%%%%%%%%%%%%%%%%%%%%%%%%%%%%%%%%%%%%%%%%%%%%%%%%%%

\mysection{Global solutions}
%%%%%%%%%%%%%%%%%%INTEGRAL INEQUALITIES%%%%%%%%%%%%%%%%%%%%%%%%%%%%%%%%%%%%%%%%%%%%%%%%%%%%%%%%%%%%%%%%%%%%%%%%%%%%%%%%%%%%%%%%%%%%%%%%%%%%%%%
\subsection{Radial solutions. Proof of Theorem D}

%%%%%%%%%%%%%%%%%%%%%%%%%%%%%%%%%%%%%%%%%%%%%%%%%%%%%%%%%%%%%%%%%%
%%%%%%%%%%%%%%%%%%%%%%%%%%%%%%%%%%%%%%%%%%%%%%%%%%%%%%%%%%%%%%%%%%

We give here the proof of Theorem D which characterizes all the positive global radial solutions of $(\ref{I-0})$ in $\BBR^N$, although Theorem C is proved in next section,  since the two proofs are completely independent. Up to translation, we  assume that the solutions are radially symmetric with respect to $0$. As for the constant $K$ in formula $(\ref{I-B-3''})$, it is given by the expression,
  \bel{2r-y}
K(N,q)=\myfrac{(1-q)\left(N-2\right)^{q-1}}{N-(N-1)q}.
\ee

%%%%%%%%%%%%%%%%%%%%%%%%%%%%%%%%%%%%%%%%%%%%%%%%%%%%%%%%%%%%%%%%%%
%%%%%%%%%%%%%%%%%%%%%%%%%%%%%%%%%%%%%%%%%%%%%%%%%%%%%%%%%%%%%%%%%%
\nind The radial form of  $(\ref{I-0})$ is the following
  \bel{2r-0}
-u''-\myfrac{N-1}{r}u'=u^p\abs{u'}^q,
\ee
and $u'(0)=0$ since any solution is $C^2$. Thus $u$ can be written under the form
  \bel{2r-1}
u(r)=u(0)+\myint{0}{r}s^{1-N}\myint{0}{s}u^p(t)\abs{u'(t)}^qt^{N-1} dt\qquad\forall r>0.
\ee
If $q\geq 1$, the solution satisfying $u(0)=a>0$ is the unique fixed point of the mapping $v\mapsto\CT[v]$ defined in the set of functions in $C([0,r_0])$ with value $a$ for $r=0$ by
$$\CT[v](r):= a+\myint{0}{r}s^{1-N}\myint{0}{s}v^p(t)\abs{v'(t)}^qt^{N-1} dt\qquad\forall r>0.
$$
Clearly $\CT$ is a strict contraction if $r_0>0$ is small enough. Since $u\equiv a$ is a solution in $\BBR^N$ it is the unique one.\smallskip

Hereafter we assume $0\leq q<1$. It can be verified that we can write $(\ref{2r-0})$ under the form
  \bel{2r-3}
\Gd_m^\gn u+(1-q)u^p=0,
\ee
where $\Gd_m^\gn u=r^{1-\gn}\left(r^{\gn-1}(u^m)'\right)'$ is the $m$-Laplacian in dimension $\gn$ applied to the radial function $u$ with $m=2-q$ and $\gn=N-(N-1)q$. An important critical value of $p$ is the following
  \bel{2r-4}
p_{crit}:=\myfrac{\gn(m-1)+m}{\gn-m}=\myfrac{\left(N-(N-1)q\right)(1-q)+2-q}{(N-2)(1-q)}.
\ee
For this specific value there exists an explicit family of ground states given by
%  \bel{2r-5}
%u_c(r)=c\left(K_c^2+r^{\frac{m}{m-1}}\right)^{\frac{m-\gn}{m+\gn}}=c\left(K_c^2+r^{\frac{2-q}{1-q}}\right)^{-\frac{(N-2)(1-q)}{2-q}},
%\ee
%where $c>0$ and
 \bel{2r-5'}
u_c(r)=c\left(Kc^{\frac{m^2}{\gn-m}}+r^{\frac{m}{m-1}}\right)^{\frac{m-\gn}{m}}=c\left(Kc^{\frac{(2-q)^2}{(N-2)(1-q)}}+r^{\frac{2-q}{1-q}}\right)^{\frac{(N-2)(q-1)}{2-q}},
\ee
with $c>0$ and
  \bel{2r-6}\BA{lll}
K=\myfrac{1}{\gn}\left(\myfrac{\gn-m}{m-1}\right)^{1-m}(1-q).
\EA\ee
Set
  \bel{2r-7}\BA{lll}
F_u(r)=r^\gn\left(\myfrac{\abs{u'}^m}{m'}+\myfrac{u^{p+1}}{p+1}+\myfrac{(\gn-m)}{m}\myfrac{\abs{u'}^{m-2}u'}{r}\right)\\[4mm]
\phantom{F_u(r)}
=r^{N-(N-1)q}\left(\myfrac{(1-q)\abs{u'}^{2-q}}{2-q}+\myfrac{u^{p+1}}{p+1}+\myfrac{(N-2)(1-q)}{2-q}\myfrac{\abs{u'}^{-q}u'}{r}\right).
\EA\ee
Then
  \bel{2r-8}\BA{lll}
F'_u(r)=r^{\gn-1}\left(\myfrac{\gn}{p+1}-\myfrac{\gn-m}{m}\right)u^{p+1}.
\EA\ee
We notice that $F'_u\equiv 0$ if and only if $p=p_{crit}$, $F'_u>0$ (resp. $F'_u<0$) if and only if $p>p_{crit}$ (resp. $p<p_{crit}$).
By \cite[Th 5.2, 5.3]{BV1}, if $p<p_{crit}$ all the solutions of $(\ref{2r-3})$ which have a finite limit at $r=0$ oscillate around $0$ when  $r\to \infty$. Hence there exists no ground state.  By \cite[Th 5.1]{BV1}, if $p>p_{crit}$, for any $\ga>0$ there exists a positive solution $u$ of $(\ref{2r-3})$ in $\BBR^N$ satisfying $u(0)=\ga$.\qeda
%%%%%%%%%%%%%%%%%%%%%%%%%%%%%%%%%%%%%%%%%%%%%%%%%%%%%%%%%%%%%%%%%%%%%%%%%%%%%%%%%%%%%%%%%%%%%%%%%%%%%%%%%%%%%%%%%%%%%%%%%%%%%%%%%%%%%%%%%%%%%%%%%%%%%%%%%%%%%%%%%%%%%%%%%%%%%%%%%%%%%%%%%%%%%%%%%%%%%%%%%
\subsection{Proof of Theorem C}
\nind{\it Step 1: Integral inequalities.} The aim of this paragraph is to prove that under the assumptions $(\ref{I-B-1})$ the gradient of any nonnegative solution $u$ of $(\ref{I-0})$ in whole $\BBR^N$ is null. The method is an extension of the one developed in \cite{GS}, \cite{BV-V}, in the sense that we still set $u=v^{-\gb}$ and $v$ satisfies $(\ref{2-1})$, and $z=\abs{\nabla v}^2$. The main novelty is that we multiply the equation satisfied by $z$ by $v^\gl z^e$ where $e>0$ and $\gl$  are two real parameters (in \cite{GS} and \cite{BV-V} they have chosen $e=0$).  The algebraic computation is heavy and we present a very technical part of it in the Appendix. Furthermore, since the exponent $e$ will sometimes take values smaller than $1$, we have to replace $z^e$ by $f(z)$ where $f$ is a smooth approximation and to consider many equations in the weak sense since $u$, and hence $v$ is merely $C^{2,q}$ if $0<q<1$. We start with the following Weintzenb\"ock inequality already used in the proof of Theorem B, but taken here in the weak sense,
$$\myint{\BBR^N}{}\left(\myfrac {1}{2}\langle\nabla z,\nabla\phi\rangle +\myfrac{(\Gd v)^2}{N}\phi
-\Gd v(\langle\nabla v,\nabla \phi\rangle+\phi\Gd v)\right)dx\leq 0,
$$
for all $\phi\in C^1_0(\BBR^N)$, $\phi\geq 0$, hence
  \bel{3-n1}\BA {lll}
\myint{\BBR^N}{}\left(\myfrac {1}{2}\langle\nabla z,\nabla\phi\rangle -\myfrac{N-1}{N}\phi(\Gd v)^2
-\Gd v\langle\nabla v,\nabla \phi\rangle\right)dx\leq 0.
\EA\ee
We choose $\phi=v^\gl f(z)\eta$ where $\eta\in C^3_0(\BBR^N)$, $\eta\geq 0$ and $f\in C^1([0,\infty))$, $f\geq 0$ and get
  \bel{3-n2}\BA {lll}
\myfrac{\gl}{2}\myint{\BBR^N}{}v^{\gl-1}f(z)\eta\langle\nabla v,\nabla z\rangle dx+\myfrac{1}{2}\myint{\BBR^N}{}v^{\gl}f'(z)\abs{\nabla z}^2\eta dx
-\myfrac{N-1}{N}\myint{\BBR^N}{}(\Gd v)^2v^\gl f(z)\eta dx\\[4mm]
\phantom{---------}
-\gl\myint{\BBR^N}{}v^{\gl-1}(\Gd v) f(z)z\eta dx-\myint{\BBR^N}{}v^{\gl}(\Gd v)f'(z)\eta\langle\nabla v,\nabla z\rangle dx\\[4mm]
\phantom{---------}
\leq -\myfrac{1}{2}\myint{\BBR^N}{}v^\gl f(z)\langle\nabla z,\nabla\eta\rangle dx
+\myint{\BBR^N}{}v^{\gl}(\Gd v)f(z)\eta\langle\nabla v,\nabla \eta\rangle dx.
\EA\ee
This inequality proved with regular functions $v$ and $f$ is extendable by density to $v\in C^2(\BBR^N)$ and $f$ locally Lipschitz continuous.
We apply this relation to the function $v$ which satisfies $(\ref{2-1})$ in the range $0\leq q<2$, $\gb\in\BBR\setminus\{1\}$ and where $s=1-q-\gb(p+q-1)$. We consider the different terms appearing in $(\ref{3-n2})$ with the help of $(\ref{2-1})$.
$$\Gs=-\myint{\BBR^N}{}(\Gd v)v^\gl f'(z)\eta\langle\nabla v\nabla z\rangle dx
=(1+\gb)\Gs_1+\abs\gb^{q-2}\gb \Gs_2
$$
with
$$\Gs_1=-\myint{\BBR^N}{}zf'(z)v^{\gl-1}\eta\langle\nabla v,\nabla\eta\rangle dx\;\;\,\text{ and }\;\;\,
\Gs_2=-\myint{\BBR^N}{}z^{\frac q2}f'(z)v^{\gl+s}\eta\langle\nabla v,\nabla\eta\rangle dx.
$$
Then
$$\BA {lll}
\Gs_1=-\myint{\BBR^N}{}\eta\langle v^{\gl-1}\nabla v,zf'(z)\nabla z\rangle dx\\[4mm]
\phantom{T_1}
=-\myint{\BBR^N}{}\eta\langle v^{\gl-1}\nabla v,\nabla(g(z)) dx\quad\text{where }\;\;g(t)=\myint{0}{t}sf'(s)ds
\\[4mm]
\phantom{\Gs_1}
=-\myint{\BBR^N}{}\langle v^{\gl-1}\nabla v,\nabla(\eta g(z))\rangle dx+\myint{\BBR^N}{}v^{\gl-1}g(z)\langle \nabla v,\nabla\eta\rangle dx
\\[4mm]
\phantom{\Gs_1}
=\myint{\BBR^N}{}\eta g(z)\nabla.(v^{\gl-1}\nabla v)dx+\myint{\BBR^N}{}v^{\gl-1}g(z)\langle \nabla v,\nabla\eta\rangle dx\\[4mm]
\phantom{\Gs_1}
=(\gl-1)\myint{\BBR^N}{}v^{\gl-2}\eta g(z)zdx+\myint{\BBR^N}{}\eta g(z)v^{\gl-1}(\Gd v) dx+\myint{\BBR^N}{}v^{\gl-1}g(z)\langle \nabla v,\nabla\eta\rangle dx\\[4mm]
\phantom{\Gs_1}
=(\gl-1)\myint{\BBR^N}{}v^{\gl-2}\eta g(z)zdx+(1+\gb))\myint{\BBR^N}{}\eta g(z)zv^{\gl-2}dx+
\\[4mm]
\phantom{\Gs_1----------}
\abs{\gb}^{q-2}\gb
\myint{\BBR^N}{}\eta g(z)z^{\frac q2}v^{\gl-1+s}dx+\myint{\BBR^N}{}v^{\gl-1}g(z)\langle \nabla v,\nabla\eta\rangle dx.
\EA$$
Similarly,
$$\BA {lll}
\Gs_2=-\myint{\BBR^N}{}\langle v^{\gl+s}\nabla v,z^{\frac{q}{2}}f'(z)\nabla z\rangle\eta dx
\\[4mm]
\phantom{\Gs_2}
=-\myint{\BBR^N}{}\langle v^{\gl+s}\nabla v,\nabla h(z)\rangle\eta dx\quad\text{where }\;\;h(t)=\myint{0}{t}s^{\frac q2}f'(s)ds
\\[4mm]
\phantom{\Gs_2}
=-\myint{\BBR^N}{}\langle v^{\gl+s}\nabla v,\nabla (\eta h(z))\rangle dx+\myint{\BBR^N}{}v^{\gl+s}h(z)\langle \nabla v,\nabla \eta\rangle  dx
\\[4mm]
\phantom{\Gs_2}
=(\gl+s)\myint{\BBR^N}{}v^{\gl+s-1}zh(z)\eta dx+\myint{\BBR^N}{}v^{\gl+s}(\Gd v)h(z)\eta dx+\myint{\BBR^N}{}v^{\gl+s}h(z)\langle \nabla v,\nabla \eta\rangle  dx\\[4mm]
\phantom{\Gs_2}
=(\gl+s+\gb+1)\myint{\BBR^N}{}v^{\gl+s-1}zh(z)\eta dx+\abs\gb^{q-2}\gb\myint{\BBR^N}{}v^{\gl+2s}z^{\frac q2}h(z)\eta dx
\\[4mm]
\phantom{\Gs_2------------------------}
+\myint{\BBR^N}{}v^{\gl+s}h(z)\langle \nabla v,\nabla \eta\rangle  dx.
\EA$$
Next we compute the term
$$\BA{lll}
\Gth_1=\myfrac{\gl}{2}\myint{\BBR^N}{}v^{\gl-1}f(z)\eta\langle\nabla v,\nabla z\rangle dx
\\[4mm]\phantom{\Gth_1}
=\myfrac{\gl}{2}\myint{\BBR^N}{}v^{\gl-1}\eta\langle\nabla v,\nabla j(z)\rangle dx\quad\text{where }\;\;j(t)=\myint{0}{t}f(s)ds
\\[4mm]\phantom{\Gth_1}
=\myfrac{\gl}{2}\myint{\BBR^N}{}v^{\gl-1}\langle\nabla v,\nabla ( j(z)\eta)\rangle dx-\myfrac{\gl}{2}\myint{\BBR^N}{}v^{\gl-1}j(z)\langle\nabla v,\nabla \eta \rangle dx
\\[4mm]\phantom{\Gth_1}
=-\myfrac{\gl}{2}\myint{\BBR^N}{}j(z)\eta\left(v^{\gl-1}\Gd v+(\gl-1)v^{\gl-2}z\right) dx-\myfrac{\gl}{2}\myint{\BBR^N}{}v^{\gl-1}j(z)\langle\nabla v,\nabla \eta \rangle dx.
\EA$$
Finally we compute
$$\BA{lll}
\Gth_2=\myfrac{1}{2}\myint{\BBR^N}{}v^{\gl}f(z)\langle\nabla z,\nabla \eta\rangle dx=\myfrac{1}{2}\myint{\BBR^N}{}v^{\gl}\langle\nabla (j(z),\nabla \eta\rangle dx
\\[4mm]\phantom{\Gth_2}
=-\myfrac{\gl}{2}\myint{\BBR^N}{}v^{\gl-1}j(z)\langle\nabla v,\nabla \eta\rangle dx
-\myfrac{1}{2}\myint{\BBR^N}{}v^{\gl}j(z)(\Gd \eta) dx.
\EA$$
Carrying forward these estimates into $(\ref{3-n2})$ yields
  \bel{3-n3}\BA {lll}
  -\myfrac{\gl(\gl-1)}{2}\myint{\BBR^N}{}v^{\gl-2}j(z)z\eta dx-\myfrac{\gl}{2}\myint{\BBR^N}{}(\Gd v)v^{\gl-1}j(z)\eta dx+\myfrac{1}{2}\myint{\BBR^N}{}v^\gl f'(z)\abs {\nabla z}^2\eta dx\\
  [4mm]\phantom{---}
  -\myfrac{N-1}{N}\myint{\BBR^N}{}(\Gd v)^2v^\gl f(z)\eta dx-\gl\myint{\BBR^N}{}(\Gd v)v^{\gl-1} zf(z)\eta dx
  \\
  [4mm]\phantom{---}
  +(1+\gb)(\gl+\gb)\myint{\BBR^N}{}v^{\gl-2}zg(z)\eta dx+(1+\gb)\abs\gb^{q-2}\gb\myint{\BBR^N}{}v^{\gl-1+s}z^{\frac q2}g(z)\eta dx
  \\ [4mm]\phantom{---}
  +\abs\gb^{q-2}\gb(\gb+1+\gl+s)\myint{\BBR^N}{}v^{\gl+s-1}zh(z)\eta dx+\gb^{2(q-1)}\myint{\BBR^N}{}v^{\gl+2s}z^{\frac q2}h(z)\eta dx
    \\ [4mm]\phantom{---}
    \leq \myint{\BBR^N}{}\left(\gl j(z)+f(z)v\Gd v-(1+\gb)g(z)-\abs\gb^{q-2}\gb v^{1+s}h(z)\right)v^{\gl-1} \langle \nabla v,\nabla\eta\rangle dx  \\ [4mm]\phantom{----------------------------}
   +\myfrac{1}{2} \myint{\BBR^N}{}v^\gl j(z)\Gd\eta dx.
  \EA\ee
  Next we fix $e\geq 0$ and choose $ f(t)=f_\varepsilon (t)=\min\{t^e+(e-1)\varepsilon ^{e}, e\varepsilon ^{e-1}t\}$ if $0\leq e< 1$ and $f(t)=t^e$ if $e\geq 1$. Then $f_\varepsilon $ is $C^1$. In order to let $\varepsilon \to 0$ in $(\ref{3-n3})$, replacing $f$, $j$, $g$ and $h$ respectively by $f_\varepsilon $, $j_\varepsilon $, $g_\varepsilon $ and $h_\varepsilon $ we notice that
$$\BA {lll}
f_\varepsilon (z)\uparrow z^e\,,\;f'_\varepsilon (z)\uparrow ez^{e-1}\,,\; g_\varepsilon (z)\uparrow \myfrac{e}{1+e}z^{1+e}\,,\; h_\varepsilon (z)\uparrow \myfrac{2e}{q+2e}z^{\frac q2+e}\,{\text{ and }}\; j_\varepsilon (z)\uparrow \myfrac{1}{1+e}z^{1+e}.
\EA$$
Since $v^\gl f'_\varepsilon (z)\abs{\nabla z}^2\eta$ converges a.e. to $v^\gl z^{e-1}\abs{\nabla z}^2\eta$ we derive by monotone convergence
$$\displaystyle \lim_{\varepsilon \to 0}\myint{\BBR^N}{}v^\gl f'_\varepsilon (z)\abs{\nabla z}^2\eta dx=e\myint{\BBR^N}{}v^\gl z^{e-1}\abs{\nabla z}^2\eta dx,
$$
which may not be finite. All the other terms in $(\ref{3-n3})$ converge by Lebesgue's theorem, hence
  \bel{3-n4}\BA {lll}\!
  -\myfrac{\gl(\gl-1)}{2(1+e)}\myint{\BBR^N}{}v^{\gl-2}z^{2+e}\eta dx-\gl\left(\myfrac{1}{2(1+e)}+1\right)\myint{\BBR^N}{}(\Gd v)v^{\gl-1}z^{1+e}\eta dx\\[4mm]\phantom{}
  \!\!\!+\myfrac{e}{2}\myint{\BBR^N}{}\!\!v^\gl z^{e-1}\abs{\nabla z}^2\eta dx-\myfrac{N-1}{N}\myint{\BBR^N}{}\!\!(\Gd v)^2v^{\gl}z^{e}\eta dx
  +\myfrac{e(1+\gb)(\gl+\gb)}{1+e}\myint{\BBR^N}{}\!\!v^{\gl-2}z^{e+2}\eta dx
  \\[4mm]\phantom{-------}
  +\abs\gb^{q-2}\gb\left(\myfrac{(\gb+1)e}{1+e}+\myfrac{2(\gb+1+\gl+s)e}{q+2e}\right)\myint{\BBR^N}{}v^{\gl+s-1}z^{\frac q2+1+e}\eta dx
    \\[4mm]\phantom{---------------------}
    +\myfrac{2e\gb^{2(q-1)}}{q+2e}\myint{\BBR^N}{}v^{\gl+2s}z^{q+e}\eta dx\leq M_0,
   \EA\ee
   where
     \bel{3-n5}\BA {lll}
     M_0=\myfrac{\gl+1+\gb}{1+e}\myint{\BBR^N}{} v^{\gl-1}z^{1+e}\langle\nabla v,\nabla\eta\rangle dx +
     \myfrac{1}{2(1+e)}\myint{\BBR^N}{} v^\gl z^{1+e}\Gd\eta dx\\[4mm]
     \phantom{M_0=-}
     +\myint{\BBR^N}{} (\Gd v)v^\gl z^{e} \langle\nabla v,\nabla\eta\rangle dx
     +\myfrac{q\abs\gb^{q-2}\gb}{q+2e}\myint{\BBR^N}{}v^{\gl+s}z\langle\nabla v,\nabla\eta\rangle dx
     \\[4mm]
     \phantom{M_0}
     \leq \myfrac{1}{2(1+e)}M+\myfrac{|\gl+1+\gb|}{1+e}L+\myfrac{|\gb|^{q-1}q}{q+2e}R,
     \EA\ee
     with
          \bel{3-n6}\BA {lll}
          M=\myint{\BBR^N}{}v^\gl z^{1+e}\abs{\Gd\eta}dx\,,\quad L=\myint{\BBR^N}{}v^{\gl-1} z^{\frac 32+e}\abs{\nabla\eta}dx\\[4mm]
       R=\myint{\BBR^N}{}v^{\gl+s} z^{\frac {q+1+2e}{2}}\abs{\nabla\eta}dx;
               \EA\ee
Here we  have used $\abs{\langle\nabla v,\nabla\eta\rangle}\leq z^\frac 12\abs{\nabla\eta}$ and  the value of $\Gd v$ given by $(\ref{2-1})$.

Set
           \bel{3-n7}\BA {lll}
G=\myint{\BBR^N}{}v^{\gl} z^{-1+e}\abs{\nabla z}^2\eta dx\,,\quad P=\myint{\BBR^N}{}v^{\gl-1+s} z^{1+\frac q2+e}\eta dx\\[4mm]
F=\myint{\BBR^N}{}v^{\gl-2} z^{2+e}\eta dx \,,\quad U=\myint{\BBR^N}{}v^{\gl+2s} z^{q+e}\eta dx.
               \EA\ee
We replace again $\Gd v$ and $(\Gd v)^2$ in the left-hand side of $(\ref{3-n4})$  to obtain
$$\myint{\BBR^N}{}(\Gd v)v^{\gl-1}z^{1+e}\eta dx=(1+\gb)F+\abs\gb^{q-2}\gb P,
$$
$$\myint{\BBR^N}{}(\Gd v)^2v^{\gl}z^{e}\eta dx=(1+\gb)^2F+\gb^{2(q-1)}U+2\abs\gb^{q-2}(1+\gb) P.
$$
 and replacing these terms into the left-hand side of $(\ref{3-n4})$, we get
          \bel{3-n8}\BA {lll}
  \myfrac{e}{2}G+A_0F+\abs{\gb}^{q} B_0P+CU\leq  \myfrac{1}{2(1+e)}M+\myfrac{|\gl+1+\gb|}{1+e}L+\myfrac{|\gb|^{q-1}q}{q+2e}R
                 \EA\ee
with%\marginpar{{\color{red} I changed $A$ to $A_0$ and $B$ to $B_0$ because $A,\ B$ are used as parameters in Step 3 and in the line previous to $(\ref{3-X})$ they are called $_0$}}
$$\BA{lll}
A_0=-\myfrac{\gl(\gl+\gb)}{2(1+e)}-\gl(1+\gb)-\myfrac{(N-1)(1+\gb)^2}{N}+\myfrac{e(1+\gb)(\gl+\gb)}{1+e}\\[4mm]
\phantom{A_0}
=-\myfrac{\gl(\gl+\gb)}{2(1+e)}-(1+\gb)\left(1-\myfrac{1+\gb}{N}+\myfrac{\gl+\gb}{1+e}\right),
\EA$$
$$\BA{lll}
\gb B_0=-\gl\left(1+\myfrac{1}{2(1+e)}\right)-\myfrac{2(N-1)(1+\gb)}{N}+\myfrac{e(1+\gb)}{1+e}+\myfrac{2e(\gb+1+\gl+s)}{q+2e}.
\EA$$
and
$$C=\gb^{2(q-1)}\left(\myfrac{2e}{q+2e}-\myfrac{N-1}{N}\right).
$$
Next we take $e=\frac{(N-1)q}{2}$ so that $C=0$ and
$$\gb B_0=-\gl\left(\myfrac {1}{N}+\myfrac{1}{2(1+e)}\right)+\myfrac{e(1+\gb)}{1+e}-\myfrac{(N-1)(\gb p+(\gb+1)q)}{N}
$$
by replacing $s$ by $1-q-\gb(p+q-1)$. Next we introduce $\gd=-\frac\gl\gb$ and $y=\frac{1+\gb}{\gb}$, hence
$$\BA {lll}
B_0=\gd\left(\myfrac{1}{N}+\myfrac{1}{2(1+e)}\right)+y\left(\myfrac{e}{1+e}-\myfrac{q(N-1)}{N}\right)-\myfrac{p(N-1)}{N}.
\EA$$
%%%%%%%%%%%%%%%%%%%%%%%%%%%%%%%%%%%%%%%%%%%%%%%%%%%%%%%%%%%%%%%%%%%%%%%%%%%%%%%%%%%%%%%%%%%%%%%%%%%%%%%%%%%%%%%%%%%%%%%%%%%%%%%%%%%%%%%%%%%%%%%%%%%%%%%%%%%%%%%%%%%%%%%%%%%%%%%%%%%%%%%%%%%%%%%%%%%%%%%%%%%%%%%%%%%%%%%%%%%%%%%%%%%%%%%%%%%%%%%%%%%%%%%%%%%%%%%%%%%%%%%%%%%%%%%%%%%%%%%%%%%%%%%%%%%%%%%%%%%%%%%%%%%%%%%%%%%%%%%%%%%%%%%%%%%%%%%%%%%%%%%%%%%%%%%%%%%%%%%%%%%%%%%%%%%%%%%%%%%%%%%%%%%%%%%%%%%%%%%%%%%%%%%%%%%%%%%%%%%%%%%%%%%%%%%%%%%%%%%%%%%%%%%%%%%%%%%%%%%%%%%%%%%%%%%%%%%%%%%%%%%%%%%%%%%%%%%%%%%%%%%%%%%%%%%%%%%%%%%%%%%%%%%%%%%%%%%%%%%%%%%%%%%%%%%%%%%%%%%%%%%%%%%%%%%%%%%%
\smallskip

%%%%%%%%%%%%%%%%%%ANALYSIS OF A_0 AND B_0%%%%%%%%%%%%%%%%%%%%%%%%%%%%%%%%%%%%%%%%%%%%%%%%%%%%%%%%%%%%%%%%%%%%%%%%%%%%%%%%%%%%%%%%%%%%%%%%%%%%%%%

\nind{\it Step 2: Study of the coefficients $A_0$ and $B_0$.}  Our method is to choose the real parameters $\gd$ and $y$ in order to ensure $A_0$ and $B_0$ to be positive. We set
\bel{3-8}\BA {lll}
m=\gd-\left(2+(N-1)q\right)y.
\EA\ee
In the sequel we keep the parameters $y$ and $m$ as variables and eliminate $\gd$.\smallskip

\nind (i) {\it Condition $A_0>0$.} We define
$$\BA {lll}
\CE(m,y)=-(2+(N-1)q)A_0\\[2mm]\phantom{--}
=\myfrac{(N-1)(Nq+1)(2+(N-1)q)}{N}y^2+2\left(1+(N-1)q\right)(m-1)y+m(m-1).
\EA$$
In the $(m,y)$-plane the set of points such that $\CE(m,y)=0$ is a conic. Its points at infinity in the associated projective space $\BBP^3(\BBR)=\{(\tilde y,\tilde m,\tilde t)\}$ satisfy, with $\tilde t=0$,
$$(2+(N-1)q)\left((N-1)q+\myfrac{N-1}{N}\right)\tilde y^2+2\left(1+(N-1)q\right)\tilde m\tilde y+\tilde m^2=0.
$$
The discriminant of this quadratic form is
$$\tilde\Gd=(1+(N-1)q)^2-(2+(N-1)q)\left((N-1)q+\myfrac{N-1}{N}\right),
$$
which is always negative since $N\geq 2$. Hence $\CE(m,y)=0$ is the equation of an ellipse, and it is easy to check that
\bel{3-9}\BA {lll}
\!\!\!\!\!\!\left\{(m,y)\in\BBR^2\!:\CE(m,y)<0\right\}\!\subset\!\left\{\!(m,y)\in\BBR^2:-\myfrac{N\left(1+(N-1)q\right)^2}{N-2+(N-1)q}<m<1\right\}\!\!.
\EA\ee
\smallskip

\nind (ii) {\it Condition $B_0>0$.} We have
$$\BA {lll}
\CD_p(m,y):=N\left(2+(N-1)q\right) B\\[1mm]\phantom{\CD_p(m,y):}
=\left(N+2+(N-1)q\right)m+2\left(N+2+(N-1)^2q\right)y-(N-1)\left(2+(N-1)q\right)p
\EA$$
The condition $\CD_p(m,y)>0$ means that $(m,y)$ belong to the upper half plane defined by the line $\CD_p$ with equation
\bel{3-10}\BA {lll}
y=-am+bp
\EA\ee
where
\bel{3-10-1}\BA {lll}
2a=\myfrac{N+2+(N-1)q}{N+2+(N^2-1)q}\;\text{ and }\;\; 2b=\myfrac{(N-1)(2+(N-1)q)}{N+2+(N^2-1)q}.
\EA\ee

The problem is reduced to find the variable $m$ so that the set $\CE(m,y)< 0$ intersects the set $y+am-bp>0$.
This means that the second degree equation $\CE(m,-am+bp)=0$ has two real zeroes. Hence
\bel{3-11}\BA {lll}
0=\myfrac{(N-1)(Nq+1)(2+(N-1)q)}{N}(am-bp)^2\\[4mm]\phantom{----------}
-2(1+(N-1)q)(m-1)(am-bp)+m^2-m
\EA\ee
where $a$ and $b$ are given by $(\ref{3-10-1})$. Its discriminant is given by $\CD=-\frac{b^2}{N}G(p,q)$ where $G(p,q)$ is defined in
$(\ref{I-B-1})$. The condition reads $G(p,q)<0$. If $q=0$ we obtain $p<\frac{N+2}{N-2}$, which the optimal condition obtained in \cite{GS}. More generally the condition on $p$ is
\bel{3-13}\BA {lll}
p<p_c(q):=\myfrac{-b(p)+\sqrt{b^2(q)+4Nq^2\left((N-1)^2q+N-2\right)}}{2\left((N-1)^2q+N-2\right)}.
\EA\ee
\smallskip

\nind{\it Step 3: Elimination of the right-hand side.}  Since $e$ can take values smaller than 1, in order to estimate $U$ in $(\ref{3-n7})$ we set $\gg_\varepsilon (z)=\min\{z^{\frac{q}{2}+e},\varepsilon ^{\frac{q}{2}+e-1}z\}$ if  $\frac{q}{2}+e<1$ and $\gg_\varepsilon (z)=z^{\frac{q}{2}+e}$ if $\frac{q}{2}+e\geq1$. We  multiply $(\ref{2-1})$ by $v^{\gl+s}\gg_\varepsilon (z)\eta$ and get
$$\BA {lll}
\abs\gb^{q-1}\myint{\BBR^N}{}v^{2s+\gl}z^{\frac q2}\gg_\varepsilon (z)\eta dx\leq\abs{\gl+s-\gb-1}\myint{\BBR^N}{}v^{s+\gl-1}z\gg_\varepsilon (z)\eta dx
+\left|\myint{\BBR^N}{}\eta v^{\gl+s}\gg'_\varepsilon (z)\langle\nabla v,\nabla z\rangle \right|dx\\[4mm]
\phantom{\abs\gb^{q-2}\gb\myint{\BBR^N}{}v^{2s+\gl}z^{\frac q2}\gg_\varepsilon (z)\eta dx}
+\left|\myint{\BBR^N}{}v^{\gl+s}\gg_\varepsilon (z)\langle\nabla v,\nabla \eta\rangle dx\right|.
\EA$$
There holds by dominated convergence
$$\displaystyle\lim_{\varepsilon \to 0}\left|\myint{\BBR^N}{}v^{\gl+s}\gg_\varepsilon (z)\langle\nabla v,\nabla \eta\rangle dx\right|
\leq \myint{\BBR^N}{}v^{\gl+s}z^{\frac{q+1}{2}+1}\abs{\nabla \eta} dx=R.
$$
We recall that $f_\varepsilon (z)=z^e$ if $e\geq 1$ and $f_\varepsilon (z)=\min\{t^e+(e-1)\varepsilon ^{e}, e\varepsilon ^{e-1}t\}$ if  $e<1$. By Cauchy-Schwarz inequality
$$\BA {lll}
\left|\myint{\BBR^N}{}\eta v^{\gl+s}\gg'_\varepsilon (z)\langle\nabla v,\nabla z\rangle \right|dx\leq
\myint{\BBR^N}{} v^{\frac\gl 2}\abs{\nabla z}\sqrt{f'_\varepsilon (z)}\myfrac{v^{\frac\gl 2+s}\gg'_\varepsilon (z)\sqrt z}{\sqrt{f'_\varepsilon (z)}}\eta dx\\[4mm]
\phantom{\myint{\BBR^N}{}\eta v^{\gl+s}\gg'_\varepsilon (z)\langle\nabla v,\nabla z\rangle dx}
\leq \ga\myint{\BBR^N}{}\myfrac{v^{\gl +2s}(\gg'_\varepsilon (z))^2z}{f'_\varepsilon (z)}\eta dx+
\myfrac{4}{\ga}\myint{\BBR^N}{}v^{\gl }\abs{\nabla z}^2f'_\varepsilon (z)\eta dx.
\EA$$
We have already seen that there exists
$$\displaystyle\lim_{\varepsilon \to 0}\myint{\BBR^N}{}v^{\gl }\abs{\nabla z}^2f'_\varepsilon (z)\eta dx=e\myint{\BBR^N}{}v^{\gl }\abs{\nabla z}^2z^{e-1}\eta dx
=eG,
$$
and $G$ is finite because of  the bounds on the right-hand side in $(\ref{3-n8})$. Considering separately the cases $0<e<\frac q2+e< 1$, $0<e<1\leq \frac q2+e$ and $1\leq e<\frac q2+e$,
we obtain, after some computations,
$$\displaystyle\lim_{\varepsilon \to 0}\myint{\BBR^N}{}\myfrac{v^{\gl +2s}(\gg'_\varepsilon (z))^2z}{f'_\varepsilon (z)}\eta dx
=\myfrac{(q+2e)^2}{4e}\myint{\BBR^N}{}v^{\gl +2s}z^{q+e}\eta dx=\myfrac{(q+2e)^2}{4e}U,
$$
This yields
\bel{3-13n}\BA {lll}
\abs\gb^{q-1} U\leq \abs{\gl+s-\gb-1}P+R +\myfrac{\ga(q+2e)^2}{4e}U+\myfrac{4e}{\ga}G.
\EA\ee
Choosing $\ga>0$ small enough we infer
\bel{3-14n}\BA {lll}
U\leq c(P+R+G),
\EA\ee
for some $c>0$ depending on the parameters.\smallskip
 %%%%%%%%%%%%%%%%%%%%%%%%%%%%%%%%%%%%%%%%%%%%%%%%%%%%%%%%%%%%%%%%%%%%%%%%%%%%%%%%%%%%%%%%%%%%%%%%%%%%%%%%%%%%%%%%%%%%%%%%%%%%%%%%%%%%%%%%%%%%%%%%%%%%%%%%%%%%%%%%%%%%%%%%%%%%%%%%%%%%%%%%%%%%%%%%%%%%%%%%%%%%%%%%%%%%%%%%%%%%%%%%%%%%%%%%%%%%%%%%%%%%%%%%%%%%%%%%%%%%%%%%%%%%%%%%%%%%%%%%%%%%%%%%%%%%%%%%%%%%%%%%%%%%%%%%%%%%%%%%%%%%%%%%%%%%%%%%%%%%%%%%%%%%%%%%%%%%%%%%%%%%%%%%%%%%%%%%%%%%%%%%%%%%%%%%%%%%%%%%%%%%%%%%%%%%%%%%%%%%%%%%%%%%%%%%%%%%%%%%%%%%%%%%

\nind { From now we  assume that the conditions on $N$, $p$ and $q$ which ensure the positivity of $A_0$ and $B_0$ are fulfilled, and that in this range of values we can find $m$ such that}
\bel{3-X}\BA {lll}
-2-(N-1)q<m<0
\EA\ee
and
\bel{3-Xm}\BA {lll}
\myfrac{{(2-q)m+2\left((2+e)p+q+e)\right)}}{p+q-1}>N.
\EA\ee

Combining $(\ref{3-14n})$ with $(\ref{3-n5})$ and $(\ref{3-n8})$ we derive,
\bel{3-17}\BA {lll}
G+P+F+U\leq  c(M+L+R),
\EA\ee
for some $c>0$ depending on $N$, $p$, $\gb$, $\gl$, $\gd$ and $q$.

The method is now to absorb the terms $M$, $L$ and $R$ by $F$, $P$ and $U$ by a repeated use of H\"older's inequality. Following the method developed in \cite {GS} and \cite{BV-V} it is simpler to return to the original function $u$ and the original exponents, we set  $\eta=\xi^\gk$, where
$$\gk=\myfrac{{(2-q)m+2\left((2+e)p+q+e)\right)}}{p+q-1}>N.$$
  Hence $(\ref{3-17})$ yields
\bel{3-18}\BA {lll}
G+\overline F+\overline P+\overline U\leq c_1(\overline M+\overline L+\overline R)
\EA\ee
where
\bel{3-19}\BA {lll}
(i)\qquad&\overline F=\myint{\BBR^N}{}u^{m-2}\abs{\nabla u}^{4+2e}\xi^\gk dx=\myint{\BBR^N}{}u^{m-2}\abs{\nabla u}^{4+(N-1)q}\xi^\gk dx
\\[4mm]
(ii)\qquad&\overline P=\myint{\BBR^N}{}u^{m+p-1}\abs{\nabla u}^{q+2+2e}\xi^\gk dx=\myint{\BBR^N}{}u^{m+p-1}\abs{\nabla u}^{2+Nq}\xi^\gk dx
\\[4mm]
(iii)\qquad&\overline U=\myint{\BBR^N}{}u^{m+2p}\abs{\nabla u}^{2q+2e}\xi^\gk dx=\myint{\BBR^N}{}u^{m+2p}\abs{\nabla u}^{(N+1)q}\xi^\gk dx
\EA\ee
and
\bel{3-20}\BA {lll}
(i)\qquad&\overline M=\myint{\BBR^N}{}u^m\abs{\nabla u}^{2+2e}\left(\xi\abs{\Gd\xi}+\abs{\nabla\xi}^2\right)\xi^{\gk-2}dx\\[4mm]
(ii)\qquad&\overline L=\myint{\BBR^N}{}u^{m-1}\abs{\nabla u}^{3+2e}\abs{\nabla\xi}\xi^{\gk-1}dx
\\[4mm]
(iii)\qquad&\overline R=\myint{\BBR^N}{}u^{m+p}\abs{\nabla u}^{q+1+2e}\abs{\nabla\xi}\xi^{\gk-1}dx.
\EA\ee

\nind{\it Absorption of $\overline L$}. In what follows $\varepsilon _i$ will denote small parameters to be fixed in order to absorb the different terms. Using H\"older's inequality we have
$$\BA{lll}
\overline L=\myint{\BBR^N}{}\left(\xi^{\ga}u^{m-1-A}\abs{\nabla u}^{B}\right)\left(\xi^{\gamma}u^{A}\abs{\nabla u}^{3+2e-B}\right)
\left(\xi^{\gk-1-\ga-\gamma}\abs{\nabla \xi}\right)dx
\\[4mm]\phantom{\overline L}
\leq \varepsilon _1^{\theta}\myint{\BBR^N}{}\xi^{\gth\ga}u^{\gth(m-1-A)}\abs{\nabla u}^{\gth B} dx+\varepsilon _2^{t}
\myint{\BBR^N}{}\xi^{t\gamma}u^{tA}\abs{\nabla u}^{t(3+2e-B)}dx
\\[4mm]\phantom{\overline L-------------------}
+\frac{1}{(\varepsilon _1\varepsilon _2)^{2\gs}}\myint{\BBR^N}{}\xi^{2\gs(\gk-1-\ga-\gamma)}\abs{\nabla \xi}^{2\gs}dx,
\EA$$
with
\bel{3-21'}\myfrac{1}{\gth}+\myfrac{1}{t}+\myfrac{1}{2\gs}=1.\ee
We choose the unknown exponents so that $\overline L\leq \varepsilon _1^{\theta}\overline F+\varepsilon _2^t\overline P+$ terms in $\xi$. We find
$$A=\myfrac{m+p-1}{t}\,,\; B=\myfrac{4+2e}{\gth},
$$
for the exponents of $\abs{\nabla u}$,
$$\ga=\myfrac{\gk}{\gth}\,,\; \gg=\myfrac{\gk}{t}\,\text{ and }\; \gk=1+\myfrac{\gk}{\gth}+\myfrac{\gk}{t}=2\gs,
$$
for the ones of $\xi$ and
$$\gth(m-1-A)=m-2\,\text{ and }\; t(3+2e-B)=q+2e+2,
$$
for the ones of $u$. Eliminating $A$ and $B$ leads to a linear system in $t$ and $\gth$,
\bel{3-21}\BA {lll}
\myfrac{m+p-1}{t}+\myfrac{m-2}{\gth}=m-1\\[4mm]
\myfrac{q+2+2e}{t}+\myfrac{4+2e}{\gth}=3+2e.
\EA\ee
The direct computation shows that
\bel{3-22}\BA {lll}
\myfrac{1}{t}=\myfrac{m+2e+2}{(2-q)m+2\left((2+e)p+q+e\right)}\,,\; \myfrac{1}{\gth}=\myfrac{(1-m)(q-1)+(3+2e)p}{(2-q)m+2\left((2+e)p+q+e\right)},
\EA\ee
hence
\bel{3-23'}\BA {lll}
2\gs=\myfrac{{(2-q)m+2\left((2+e)p+q+e)\right)}}{p+q-1}=\gk.
\EA\ee
We set $Y=2\gs (p+q-1))=\gth X$, hence
$$Y=(2-q)m+(N-1)pq+4p+(N+1)q
$$
and
$$X=(1-m)(q-1)+(N-1)pq+3p.
$$
First we check that under $(\ref{3-X})$ $Y>0$. Indeed we have
\bel{3-23}\BA {lll}Y-2(p+q-1)=(2-q)m+(N-1)pq+2p+(N-1)q+2\\\phantom{Y-2(p+q-1)}
\geq (2-q)m+(N-1)pq+2p-m\\\phantom{Y-2(p+q-1)}
=(1-q)m+p(2+(N-1)q)>(1-q-p)m>0,
\EA\ee
from $(\ref{3-X})$.  Next we check that $X>0$: if $q\geq 1$ this is clear by $(\ref{3-X})$. If $0\leq q<1$ we have, also by $(\ref{3-X})$,
$$X\geq -(3+(N-1)q)(1-q)+((N -1)q + 3)p=(3+(N-1)q)(p+q-1)>0.
$$
As a by-product we derive that $\gth$ and $t$ are positive and therefore larger than $1$ because of $(\ref{3-21'})$. \smallskip

\nind{\it Absorption of $\overline R$}. We introduce new parameters $A$, $B$, $t$, $\gth$ in order to absorb $\overline R$ by $\overline P+\overline U+$ term in $\xi$.
$$\BA{lll}
\overline R=\myint{\BBR^N}{}\left(\xi^{\ga}u^{m+p-A}\abs{\nabla u}^{B}\right)\left(\xi^{\gamma}u^{A}\abs{\nabla u}^{q+1+2e-B}\right)
\left(\xi^{\gk-1-\ga-\gamma}\abs{\nabla \xi}\right)dx
\\[4mm]\phantom{\overline R}
\leq \varepsilon _3^\theta\myint{\BBR^N}{}\xi^{\gth\ga}u^{\gth(m+p-A)}\abs{\nabla u}^{\gth B} dx+
\varepsilon _4^t\myint{\BBR^N}{}\xi^{t\gamma}u^{tA}\abs{\nabla u}^{t(q+1+2e-B)}dx
\\[4mm]\phantom{\overline R-------------------}
+\frac{1}{(\varepsilon _3\varepsilon _4)^{2\gs}}\myint{\BBR^N}{}\xi^{2\gs(\gk-1-\ga-\gamma)}\abs{\nabla \xi}^{2\gs}dx,
\EA$$
with $t,\gth$ and $\gs$ satisfying $(\ref{3-21'})$. Hence
\bel{3-24}\BA {lll}
\phantom{--}(m+p-A)\gth&\!\!=m+p-1\\
\phantom{(m+p^,,,,A)\gth}
B\gth&\!\!=q+2+2e\\\phantom{(m+p^,,,,A)\gth}
At&\!\!=m+2p\\\phantom{}
(q+1+2e-B)t&\!\!=2q+e.
\EA\ee
Thus
$$\BA {lll}
\myfrac{m+2p}{t}+\myfrac{m+p-1}{\gth}=m+p\\[4mm]
\myfrac{2q+e}{t}+\myfrac{q+2+2e}{\gth}=q+1+2e
\EA$$
{\it Mutatis mutandis} it yields, always with $\gk=2\gs$ given by $(\ref{3-23'})$, and now
\bel{3-25}\BA {lll}
\myfrac{1}{\gth}=\myfrac{m(1-q)+2(e+1)p}{Y}\,,\; \myfrac{1}{t}=\myfrac{m+p+q+1+2e}{Y},
\EA\ee
where $Y$ is unchanged. Condition $(\ref{3-X})$ implies $m+p+q+1+2e\geq p+q-1>0$, hence $t>0$. Furthermore
$$m(1-q)+2(e+1)p>(m+2+2e)(1-q)=(m+2+(N-1)q)(1-q)>0.
$$
which implies $\gth>0$ if $q\leq 1$. If $q>1$ and since $m<0$, then
$$m(1-q)+2(e+1)p>2(e+1)p>0,
$$
which again yields  $\gth>0$.\smallskip

\nind{\it Absorption of $\overline M$}. We set
$$\BA{lll}
\overline M_1=\myint{\BBR^N}{}u^m\abs{\nabla u}^{2+2e}\abs{\Gd\xi}\xi^{\gk-1}dx
\\[4mm]\phantom{\overline M_1}
=\myint{\BBR^N}{}\left(\xi^{\ga}u^{m-A}\abs{\nabla u}^{B}\right)\left(\xi^{\gamma}u^{A}\abs{\nabla u}^{2+2e-B}\right)
\left(\xi^{\gk-1-\ga-\gamma}\abs{\Gd \xi}\right)dx
\\[4mm]\phantom{\overline M_1}
\leq \varepsilon _5^\theta\myint{\BBR^N}{}\xi^{\gth\ga}u^{\gth(m-A)}\abs{\nabla u}^{\gth B} dx+
\varepsilon _6^t\myint{\BBR^N}{}\xi^{t\gamma}u^{tA}\abs{\nabla u}^{t(2+2e-B)}dx
\\[4mm]\phantom{\overline M_1-------------------}
+\frac{1}{(\varepsilon _5\varepsilon _6)^{\gs}}\myint{\BBR^N}{}\xi^{\gs(\gk-1-\ga-\gamma)}\abs{\Gd \xi}^{\gs}dx,
\EA$$
with
\bel{3-26}\BA {lll}
\myfrac{1}{\gth}+\myfrac{1}{t}+\myfrac{1}{\gs}=1.
\EA\ee
If we try to absorb $\overline M_1$ by $\overline F+\overline U+$ term in $\xi$, we obtain
\bel{3-27}\BA {lll}
\phantom{-B-}(m-A)\gth&\!\!=m-2\\
\phantom{(m+^,,,,A)\gth}
B\gth&\!\!=4+2e\\\phantom{(m+^,,,,A)\gth}
At&\!\!=m+2p\\\phantom{A}
(2+2e-B)t&\!\!=2q+2e.
\EA\ee
We find
\bel{3-28}\BA {lll}
\myfrac{1}{\gth}=\myfrac{m(1-q)+2(e+1)p}{Y}\,,\; \myfrac{1}{t}=\myfrac{m+2+2e}{Y}.
\EA\ee
Clearly $(\ref{3-26})$ holds and conditions $m(1-q)+2(e+1)p>0$ and $m+2+2e>0$ are satisfied under the same condition as for the treatment of $\overline R$. The same proof works for $\overline M_2=\myint{\BBR^N}{}u^m\abs{\nabla u}^{2+2e}\abs{\nabla\xi}^2\xi^{\gk-2}dx$.\medskip

\medskip

\nind{\it Step 4: End of the proof}. It follows from Step 3 by choosing the parameters $\varepsilon _i$ small enough that there holds for any nonnegative $\xi\in C^\infty_0(\BBR^N)$:
\bel{3-29}\BA {lll}
\overline L+\overline F+\overline U\leq c_3\myint{\BBR^N}{}\left(\abs{\Gd \xi}^{\gs}+\abs{\nabla \xi}^{2\gs}\right)dx.
\EA\ee
Assuming $\xi$ has support in $B_1$ and applying $(\ref{3-29})$ to $\xi_R:x\mapsto\xi(\frac{x}{R})$, we derive
\bel{3-30}\BA {lll}
\overline L+\overline F+\overline U\leq c_3R^{N-2\gs}\myint{B_1}{}\left(\abs{\Gd \xi}^{\gs}+\abs{\nabla \xi}^{2\gs}\right)dx.
\EA\ee
Hence, as $2\gs>N$ we infer $\overline L+\overline F+\overline U=0$ by letting $R\to\infty$. It remains to prove that such an estimate holds. The condition
$2\gs>N$ is equivalent to
\bel{3-31}\BA {lll}
m(2-pq+(N-1)pq+N+q>(N-4)p,
\EA\ee
or, equivalently,
\bel{3-32}\BA {lll}
(2-q)(m+2+(N-1)q)>(p+q-1)(N-4-(N-1)q).
\EA\ee
Since the left-hand side is positive by $(\ref{3-X})$, this inequality holds at least when
\bel{3-33}\BA {lll} N=3,4\,\text{ or }\, N\geq 5\, \text{ and }\, q\geq \myfrac{N-4}{N-1}.
\EA\ee
The general proof of $(\ref{3-X})$, $(\ref{3-32})$ is technical and given in Appendix.\qeda
%%%%%%%%%%%%%%%%%%%%%%%%%%%%%%%%%%%%%%%%%%%%%%%%%%%%%%%%%%%%%%%%%%%%%%%%%%%%that $X$, $P$ and $U$ are linked by the relation%%%%%%%%%%%%%%%%%%%%%%%%%%%%%%%%%%%%%%%%%%%%%%%%%%%%%%%%%%
%%%%%%%%%%%%%%%%%%%%%%%%%%%%%%%%%%%%%%%%%%%%%%%%%%%%%%%%%%%%%%%%%%%%%%%%%%%%%%%%%%%%%%%%SEPARABLE%%SOLUTIONS%%%%%%%%%%%%%%%%%%%%%%%%%%%%%%%%%%%%%%%%%%%%%%%%%%%%%%%%%%%%%%%%%%%%%%%%%%%%%%%%%%%%%%%%%%%%%%%%%%%%%%%%%%%%%%%%%%%%%%%%%%%%%%%%%%%%%%%%%%%%%%%%%%%%%%%%%%%%%%%%%%%%%%%
%%%%%%%%%%%%%%%%%%%%%%%%%%%%%%%%%%%%%%%%%%%%%%%%%%%%%%%%%%%%%%%%%%%%%%%%%%%%%%%%%%%%%%%%%%%%%%%%%%%%%%%%%%%%%%%%%%%%%%%%%%%%%%%%%%%%%%
\mysection{Separable solutions}

In the sequel we set $n=N-1$, and consider a more general equation on $S^n$,
\bel{L-1}
-\Gd'\gw+\mu\gw=\gw^p(\gg^2\gw^2+\abs{\nabla'\gw}^2)^{\frac q2},
\ee
where $\gg>0$ and $\mu$ are parameters and $\Gd'$ and $\nabla'$ are respectively the Laplace-Beltrami operator and the covariant gradient, which can be assimilated to the tangential gradient on $S^n$. Notice that if $\mu>0$
there exists a constant solution $\gw_\mu$ to $(\ref{L-1})$ given by
\bel{L-2}
\gw_\mu=\left(\frac{\mu}{\gg^q}\right)^{\frac{1}{p+q-1}}.
\ee
%%%%%%%%%%%%%%%%%%%%%%%%%%%%%%%%%%%%%%%%%%%%%%%%%%%%%%%%%%%%%%%%%%%%%%%%%%%%%%%%%%%%%%%%%%%%%%%%%%%%%%%%%%%%%%%%%%%%%%%%%%%%%%%%%%%%%%UNIFORM%BOUNDS%%%%%%%%%%%%%%%%%%%%%%%%
%%%%%%%%%%%%%%%%%%%%%%%%%%%%%%%%%%%%%%%%%%%%%%%%%%%%%%%%%%%%%%%%%%%%%%%%%%%%%%%%%%%%%%%%%%%%%%%%%%%%%%%%%%%%%%%%%%%%%%%%%%%%%%%%%%%%%%%%%%%%%%%%%%%%%%%%%%%%%%%%%%%%%%%%%%%%%%

%%%%%%%%%%%%%%%%%%%%%%%%%%%%%%%%%%%%%%%%%%%%%%%%%%%%%%%%%%%%%%%%%%%%%%BOUNDS%%ON%%S^N_+%%%%%%%%%%%%%%%%%%%%%%%%%%%%%%%%%%%%%%%%%%%%%%%%%%%%%%%%%%%%%%%%%%%%%%%%%%%%%%%%%%%%%%%%%%%%%%%%%%%%%%%%%%%%%%%%%%%%%%%%%%%%%%%%%%%%%%%%%%%%%%%%%%%%%%%%%%%%%%%%%%%%%%%%%%%%%%%%%%%%%%%%%%%%%%%%%%%%%%%%%%%%%%%%%%%%%%%%%%%%%%%%%%%%%%%%%%%%%%%%%%%%%%%%%%%%%%%%%%%%%%%%%
%%%%%%%%%%%%%%%%%%%%%%%%%%%%%%%%%%%%%%%%%%%%%%%%%%%%%%%%%%%%%%%%%%%%%%%%%%%%%%%%%%%%%%%%%%%%%%%%%%%%%%%%%%%%%%%%%%%%

%%%%%%%%%%%%%%%%%%%%%%%%%%%%%%%%%%%%%%%%%%%%%%%%%%%%%%%%%%%%%%%%%%%%%%BOUNDS%%ON%%S^N%%%%%%%%%%%%%%%%%%%%%%%%%%%%%%%%%%%%%%%%%%%%%%%%%%%%%%%%%%%%%%%%%%%%%%%%%%%%%%%%%%%%%%%%%%%%%%%%%%%%%%%%%%%%%%%%%%%%%%%%%%%%%%%%%%%%%%%%%%%%%%%%%%%%%%
\subsection{Uniform bounds: Proof of Theorem E}

We set $a\vee b=\max\{a,b\}$ and $a\wedge b=\min\{a,b\}$.
By integration on $S^n$ and H\"older's inequality there holds
$$\BA {llll}\mu\abs{S^n}^{\frac{p+q-1}{p+q}}\left(\myint{S^n}{}\gw^{p+q}dS\right)^{\frac1{p+q}}\geq \mu\myint{S^n}{}\gw dS=\myint{S^n}{}\gw^p(\gg^2\gw^2+\abs{\nabla'\gw}^2)^{\frac{q}{2}} dS\geq\gg^q\myint{S^n}{}\gw^{p+q}dS.
\EA$$
Hence
\bel{L-3}\BA {ll}
\norm\gw_{L^{p+q}}\leq \left(\myfrac{\mu}{\gg^q}\right)^{\frac1{p+q-1}}\abs{S^n}^{\frac{1}{p+q}}.
\EA\ee
Therefore

\bel{L-4}\BA {ll}
\mu\myint{S^n}{}\gw dS=\myint{S^n}{}\gw^p(\gg^2\gw^2+\abs{\nabla'\gw}^2)^{\frac{q}{2}} dS\leq \left(\myfrac{\mu^{p+q}}{\gg^q}\right)^{\frac1{p+q-1}}\abs{S^n}.
\EA\ee
For $\ga>0$, we also have
\bel{L-4-1}\BA {ll}
\myint{S^n}{}\gw^{p+\ga}\left(\gg^2\gw^2+\abs{\nabla'\gw}^2\right)^{\frac q2}dS=\myint{S^n}{}\left(\ga\gw^{\ga-1}\abs{\nabla'\gw}^2+\mu\gw^{\ga+1}\right)dS\\[4mm]
\phantom{\myint{S^n}{}\gw^{p+\ga}\left(\gg^2\gw^2+\abs{\nabla'\gw}^2\right)^{\frac q2}dS}
\geq \mu\wedge\frac{4\ga}{(\ga+1)^2}\myint{S^n}{}\!\left(\abs{\nabla'\gw^{\frac{\ga+1}{2}}}^2\!\!+(\gw^{\frac{\ga+1}{2}})^2\right)dS
\\[4mm]
\phantom{\myint{S^n}{}\gw^{p+\ga}\left(\gg^2\gw^2+\abs{\nabla'\gw}^2\right)^{\frac q2}dS}
\geq C_1\left(\mu\wedge\frac{4\ga}{(\ga+1)^2}\right)\norm\gw^{\ga+1}_{L^{\frac{n(\ga+1)}{n-2}}},
\EA\ee
using Sobolev inequality in $H^1(S^n)$. Furthermore, by H\"older's inequality,
\bel{L-4-2}\BA {ll}
\myint{S^n}{}\gw^{p+\ga}\left(\gg^2\gw^2+\abs{\nabla'\gw}^2\right)^{\frac q2}dS=
\myint{S^n}{}\gw^{p+\ga+q\frac{1-\ga}{2}}\left(\gg^2(\gw^{\frac{\ga+1}{2}})^2+\frac{4}{(\ga+1)^2}\abs{\nabla'\gw^{\frac{\ga+1}{2}}}^2\right)^{\frac q2}dS\\[4mm]
\phantom{----,--}
\leq \gg^q\vee (\frac{2}{1+\ga})^q\left(\myint{S^n}{}\gw^{\frac{2p+q}{2-q}+\ga}dS\right)^{\frac{2-q}{2}}
\left(\myint{S^n}{}\!\left(\abs{\nabla'\gw^{\frac{\ga+1}{2}}}^2\!\!+(\gw^{\frac{\ga+1}{2}})^2\right)dS\right)^{\frac{q}{2}}.
\EA\ee
It implies
\bel{L-4-3}\BA {ll}
\myint{S^n}{}\!\left(\abs{\nabla'\gw^{\frac{\ga+1}{2}}}^2\!\!+(\gw^{\frac{\ga+1}{2}})^2\right)dS
\leq \left(\gg\vee\frac{2}{1+\ga}\right)^{\frac{2q}{2-q}}\myint{S^n}{}\gw^{\frac{2p+q}{2-q}+\ga}dS.
\EA\ee
Jointly with $(\ref{L-4-1})$ it yields
\bel{L-4-4}\BA {ll}
\norm\gw^{\ga+1}_{L^{\frac{n(\ga+1)}{n-2}}}
\leq C_2\myfrac{\left(\gg\vee(1+\ga)^{-1}\right)^{\frac{2q}{2-q}}}{\mu\wedge\ga(\ga+1)^{-2}}\myint{S^n}{}\gw^{\frac{2p+q}{2-q}+\ga}dS.
\EA\ee
We define the sequence $\{\ga_k\}$ by
\bel{L-4-5}\BA {ll}
\myfrac{2p+q}{2-q}+\ga_k=\myfrac{n(\ga_{k-1}+1)}{n-2}\Longleftrightarrow \ga_k+1=\myfrac{n(\ga_{k-1}+1)}{n-2}-2\myfrac{p+q-1}{2-q}.
\EA\ee
The value of $\ga_0>0$ will be made precise later on.  The value of $\ga_k$ is explicit:
\bel{L-4-7}\BA {ll}
\ga_k+1=\left(\myfrac{n}{n-2}\right)^k\left(\ga_0+1-\myfrac{(p+q-1)(n-2)}{2-q}\right)+\myfrac{(p+q-1)(n-2)}{2-q}.
\EA\ee
Notice that since $(n-2)p+(n-1)q<n$, then  $1-\frac{(p+q-1)(n-2)}{2-q}>0$. Asymptotically
\bel{L-4-8}\BA {ll}
\ga_k+1=A\ell^k+O(1)\quad\text{with }\,\ell=\myfrac{n}{n-2}>1.
\EA\ee
We set $X_k=\norm\gw_{L^{\frac{n(\ga_k+1)}{n-2}}}$, hence $(\ref{L-4-4})$ reads
\bel{L-4-9}\BA {ll}
X_k\leq \left(C_2\myfrac{\left(\gg\vee(1+\ga_k)^{-1}\right)^{\frac{2q}{2-q}}}{\mu\wedge\ga_k(\ga_k+1)^{-2}}\right)^{\frac{1}{1+\ga_k}}
X_{k-1}^{1+2\frac{p+q-1}{(2-q)(\ga_k+1)}}.
\EA\ee
Because of $(\ref{L-4-8})$,
$$\BA{lll}\displaystyle
C_2\myfrac{\left(\gg\vee(1+\ga_k)^{-1}\right)^{\frac{2q}{2-q}}}{\mu\wedge\ga_k(\ga_k+1)^{-2}}
\leq C_3\ga_k^{-1}\gg^{\frac{2q}{2-q}}\\[4mm]\phantom{-----}\displaystyle
\Longrightarrow \left(C_2\myfrac{\left(\gg\vee(1+\ga_k)^{-1}\right)^{\frac{2q}{2-q}}}{\mu\wedge\ga_k(\ga_k+1)^{-2}}\right)^{\frac{1}{1+\ga_k}}\leq C^{\frac{1}{1+\ga_k}}_4\gg^{\frac{2q}{(1+\ga_k)(2-q)}}
\EA$$
provided $\ga_0\geq\ge_0>0$, hence
\bel{L-4-10}\BA {ll}
X_k\leq C^{\frac{1}{1+\ga_k}}_4\gg^{\frac{2q}{(1+\ga_k)(2-q)}}
X_{k-1}^{1+2\frac{p+q-1}{(2-q)(\ga_k+1)}}.
\EA\ee
 Next we construct by induction an increasing sequence $\Gg_k$ such that
 \bel{L-4-11}\BA {ll}
X_{k-1}\leq \Gg_{k-1}\gg^{-\frac{q}{p+q-1}}.
\EA\ee
Then $(\ref{L-4-11})$ holds at the order $k$ with
 \bel{L-4-12}\BA {ll}
\Gg_k= C^{\frac{1}{1+\ga_k}}_4\Gg_{k-1}^{1+2\frac{p+q-1}{(2-q)(\ga_k+1)}},
\EA\ee
and $\Gg_0$ will be fixed later on. We can assume that $C_4\geq 1$, therefore $\{\Gg_k\}$ is increasing. If we put $\gth_k=\ln\Gg_k$, then
$$\gth_k=\frac{1}{1+\ga_k}\ln C_4+\left(1+2\frac{p+q-1}{(2-q)(\ga_k+1)}\right)\gth_{k-1}= \myfrac{A}{1+\ga_k}+\left(1+\myfrac{B}{1+\ga_k}\right)\gth_{k-1}
$$
Put $\tilde\gth_k=\gth_k+\frac{A}{B}$, then
$$\tilde\gth_k=\left(1+\myfrac{B}{1+\ga_k}\right)\tilde\gth_{k-1}\Longrightarrow \tilde\gth_k=\tilde\gth_0\prod_{j=1}^k\left(1+\myfrac{B}{1+\ga_j}\right)
$$
Finally
 \bel{L-4-13}\gth_k=\left(\gth_0+\frac AB\right)\prod_{j=1}^k\left(1+\myfrac{B}{1+\ga_j}\right)-\frac AB,
\ee
 and we conclude that
 \bel{L-4-14}\Gg^*=\lim_{k\to\infty}\Gg_k.
\ee
By standard linear elliptic regularity theory with $L^1$ data and $(\ref{L-4})$,
\bel{En-5-0}\BA {ll}
\norm{\gw}_{L^{\frac{n}{n-2},\infty}}+\norm{\nabla'\gw}_{L^{\frac{n}{n-1},\infty}}\leq C_5
\myint{S^n}{}\left(\gw^p(\gg^2\gw^2+\abs{\nabla'\gw}^2)^{\frac{q}{2}} +\mu\gw\right)dS\\[4mm]
\phantom{\norm{\gw}_{L^{\frac{n}{n-2},\infty}}+\norm{\nabla'\gw}_{L^{\frac{n}{n-1},\infty}}}
\leq 2C_5\left(\myfrac{\mu^{p+q}}{\gg^q}\right)^{\frac{1}{p+q-1}}\abs{S^n},
\EA\ee
where  $L^{r,\infty}$ denotes the usual Marcinkiewicz spaces (or Lorentz spaces). For any $1<\gt<\frac{n}{n-1}$, there exists $C_6=C(n,\gt)$ such that
$$\norm\gw_{W^{1,\gt}}\leq C_6\left(\norm{\gw}_{L^{\frac{n}{n-2},\infty}}+\norm{\nabla'\gw}_{L^{\frac{n}{n-1},\infty}}\right),
$$
and by Sobolev inequality
$$\norm\gw_{L^{\gt^*}}\leq C_7\norm\gw_{W^{1,\gt}},
$$
where $\frac{1}{\gt^*}=\frac{1}{\gt}-\frac{1}{n}$. Since $n(p-2)+q(n-1)<n$ is equivalent to $\frac{2p+q}{2-q}<\frac{n}{n-2}$, we can take
$\frac{1}{\gt}=\left(\frac{2-q}{2p+q}+\ga_1\right)^{-1}+\frac 1n$ for some $\ga_1>0$. Using $(\ref{L-4-5})$ we define the initial data of $\{\ga_k\}$ by
$$\frac{2p+q}{2-q}+\ga_1=\frac{n(\ga_0+1)}{n-2}
$$
and we derive
\bel{En-5-1}\BA {ll}
X_0=\norm\gw_{L^{\frac{n(\ga_0+1)}{n-2}}}=\norm\gw_{L^{\frac{2p+q}{2-q}+\ga_1}}\leq C_8\left(\myfrac{\mu^{p+q}}{\gg^q}\right)^{\frac{1}{p+q-1}}.
\EA\ee
Finally, we can fix
$$e^{\gth_0}=\Gg_0=C_8\mu^{\frac{p+q}{p+q-1}},$$
and derive from $(\ref{L-4-13})$
\bel{En-5-2}\BA {ll}
\Gg^*=C_9e^{\gth_0\frac{A}{B}\prod_{j=0}^\infty(1+\frac{B}{1+\ga_j})}=c_{10}\mu^{\frac{p+q}{p+q-1}\frac{A}{B}\prod_{j=0}^\infty(1+\frac{B}{1+\ga_j})}
\EA\ee
with
$$C_9=e^{\frac{A}{B}\prod_{j=0}^\infty(1+\frac{B}{1+\ga_j})-\frac{A}{B}}\;\text{ and } \;\;c_{10}=C_9C^{\frac{A}{B}\prod_{j=0}^\infty(1+\frac{B}{1+\ga_j})}_8.
$$
Since
$$\norm\gw_{L^\infty}\leq \gg^{-\frac{q}{p+q-1}}\Gg^*,
$$
the estimate follows.\qeda
\medskip
%%%%%%%%%%%%%%%%%%%%%%%%%%%%%%%%%%%%%%%%%%%%%%%%%%%%%%%%%%%%%%%%%%%%%%%%%%%%%%%%%%%%%%%%%%%%%%%%%%%%%%%%%%%%%%%%%%%%%%%%%%%%%%%%%%%%%%
%%%%%%%%%%%%%%%%%%%%%%%%%%%%%%%%%%%%%%%%%%%%%%%%%%%%%%%%%%%%%%%%%%%%%%%%%%%%%%%%%%%%%%%%%%%%%%%%%%%%%%%%%%%%%%%%%%%%%%%%%%%%%%%%%%%%%%
%%%%%%%%%%%%%%%%%%%%%%%%%%%%%%%%%%%%%%%%%%%%%%%%%%%%%%%%%%%%%%%%%%%%%%%%%%%%%%%%%%%%%%%%%%%%%%%%%%%%%%%%%%%%%%%%%%%%%%%%%%%%%%%%%%%%%%
%%%%%%%%%%%%%%%%%%%%%%%%%%%%%%%%%%%%%%%%%%%%%%%%%%%%%%%%%%%%%%%%%%%%%%%%%%%%%%%%%%%%%%%%%%%%%%%%%%%%%%%%%%%%%%%%%%%%%%%%%%%%%%%%%%%%%%
%%%%%%%%%%%%%%%%%%%%%%%%%%%%%%%%%%%%%%%%%%%%%%%%%%%%%%%%%%%%%%%%%%%%%%%%%%%%%%%%%%%%%%%%%%%%%%%%%%%%%%%%%%%%%%%%%%%%%%%%%%%%%%%%%%%%%%
%%%%%%%%%%%%%%%%%%%%%%%%%%%%%%%%%%%%%%%%%%%%%%%%%%%%%%%%%%%%%%%%%%%%%%%%%%%%%%%%%%%%%%%%%%%%%%%%%%%%%%%%%%%%%%%%%%%%%%%%%%%%%%%%%%%%%%
%%%%%%%%%%%%%%%%%%%%%%%%%%%%%%%%%%%%%%%%%%%%%%%%%%%%%%%%%%%%%%%%%%%%%%%%%%%%%%%%%%%%%%%%%%%%%%%%%%%%%%%%%%%%%%%%%%%%%%%%%%%%%%%%%%%%%%

\nind\Remark We conjecture that the best exponent $a$ is equal to $1$ and
\bel{En-5-3}\norm\gw_{L^\infty}\leq C_{11}\left(\frac{\mu}{\gg^q}\right)^{\frac{1}{p+q-1}}.
\ee
Notice that there always holds
\bel{En-5-4}\BA {lll}\displaystyle
\min_{S^n}\gw\leq \left(\frac{\mu}{\gg^q}\right)^{\frac{1}{p+q-1}}\leq \max_{S^n}\gw.
\EA
\ee

%%%%%%%%%%%%%%%%%%%%%%%%%%%%%%%%%%%%%%%%%%%%%%%%%%%%%%%%%%%%%%%%%%%%RIGIDITY%%%%%%%%%%RIGIDITY AND SYMMETRY%%%%%%%%%%%%%%%%%%%%%%%%%%%%%%%%%%%%%%%%%%%%%%%%%%%%%%%%%%%%%%%%%%%%%%%%%%%%%%%%%%%%%%%%%%%%%%%%%%%%%%%%%%%%%%%%%%%%%%%%%%%%%%%%%%%%%%%%%%%%%%%%%%%%%%%%%%%%%%%%%%%%%%%%%%%%%%%%%%%%%%%%%%%%%%%%%%%%%%%%%%%%%%%%%%%%%%%%%%%%%%%%%%%%%%%%%%%%%%%%%
%%%%%%%%%%%%%%%%%%%%%%%%%%%%%%%%%%%%%%%%%%%%%%%%%%%%%%%%%%%%%%%%%%%%Symmetry%%%%%%%%%%%%%%%%%%%%%%%%%%%%%%%%%%%%%%%%%%%%%%%%%%%%%%%%%%%%%%%%%%%%%%%%%%%%%%%%%%%%%%%%%%%%%%%%%%%%%%%%%%%%%%%%%%%%%%%%%%%%%%%%%%%%%%%%%%%%%%%%%%%%%%%%%%%%%%%%%%%%%%%%%%%%%%%%%%%%%%%%%%%%%%%%%%%%%%%%%%%%%%%%%%%%%%%

\subsection{Rigidity and symmetry}

\bth{rig-1} Assume $\gg$, $\mu>0$, $p+q-1>0$ and $\gw$ is a solution of $(\ref{L-1})$ on $S^n$ such that
\bel{constr}\BA {lll}
(i)\qquad\qquad&\gg^2\gw^2+\abs{\nabla'\gw}^2\leq c_1^2&\qquad\text{if }\;p\geq 1),\qquad\qquad\qquad\\[2mm]
(ii)\qquad\qquad&
c_2^2\leq\gg^2\gw^2+\abs{\nabla'\gw}^2\leq c_1^2&\qquad\text{if }\;0\leq p< 1),\qquad\qquad\qquad
\EA\ee
for some $c_1,c_2>0$ and set
\bel{cstar}
c_*=\left\{\BA {ll} c_1\quad&\text{if }p\geq 1,
\\[2mm]
c_2^{\frac{p-1}{p+q-1}}c_1^{\frac{q}{p+q-1}}\quad&\text{if }0\leq p < 1.
\EA\right.
\ee
If $c^*$ satisfies
\bel{cstar-1}
c_*^{p+q-1}\leq \myfrac{2(n+\mu)}{q\gg^{-p}\sqrt{n}+2(p+q)\gg^{1-p}},
\ee
then $\gw$ is constant.
\es
\Proof
If $w$ is a function defined on $S^{n}$, we put
$$\bar w=\myfrac{1}{|S^{n}|}\myint{S^{n}}{}w(\gs)dS.
$$
If
\bel{R0}-\Gd'\gw+\mu\gw-|\gw|^{p-1}\gw(\gg^2\gw^2+|\nabla'\gw|^2)^{\frac{q}{2}}=0,
\ee
we have
$$-\Gd'\bar \gw +\mu\bar\gw-\overline{|\gw|^{p-1}\gw(\gg^2\gw^2+|\nabla'\gw|^2)^{\frac{q}{2}}}=0.
$$
Since $\overline{w-\overline w}=0$, $\overline w$ in the orthogonal projection in $L^2(S^n)$ of $w$ on $ker (-\Gd')$ and $n$ is the first nonzero eigenvalue, we have
\bel{R1}\BA {ll}\myint{S^{n}}{}\left(\gw^p(\gg^2\gw+|\nabla'\gw|^2)^{\frac{q}{2}}-\overline{\gw^p(\gg^2\gw^2+|\nabla'\gw|^2)^{\frac{q}{2}}} \right)(\gw-\bar\gw)dS\\[4mm]\phantom{---}
=\myint{S^{n}}{}\left(\gw^p(\gg^2\gw^2+|\nabla'\gw|^2)^{\frac{q}{2}}-\bar \gw^p(\gg^2\bar \gw^2+|\nabla'\bar \gw|^2)^{\frac{q}{2}} \right)(\gw-\bar\gw)dS\EA\ee
and
\bel{R2}\BA {lll}\myint{S^{n}}{}\left(-\Gd'(\gw-\bar\gw+\mu(\gw-\bar\gw))\right)(\gw-\bar\gw)dS
=\myint{S^{n}}{}\left(|\nabla'(\gw-\bar\gw)|^2+\mu(\gw-\bar\gw)^2)\right)dS
\\[4mm]\phantom{\myint{S^{n}}{}\left(-\Gd'(\gw-\bar\gw+\mu(\gw-\bar\gw))\right)(\gw-\bar\gw)dS}
\geq (\mu+n)\myint{S^{n}}{}(\gw-\bar\gw)^2dS
\EA
\ee
Set $F(X,Y)=\abs X^{p-1}X(\gg^2X^2+\abs{Y}^2)^{\frac q2}$ and
$$G:=\{(X,Y)\in \BBR\ti\BBR^{n}:c^2_2\leq\gg^2X^2+\abs{Y}^2\leq c^2_1\}.$$
Then
$$D_1F(X,Y)=\abs{X}^{p-1}\left(\gg^2X^2+\abs{Y}^2\right)^{\frac{q}{2}-1}\left(\gg^2(p+q)\abs{X}^{2}+p\abs{Y}^2\right)$$
$$D_2F(X,Y)=q\abs X^{p-1}X\left(\gg^2X^2+\abs{Y}^2\right)^{\frac{q}{2}-1}Y
$$
If we assume that $(\gw,\nabla'\gw)\in G$, we have, with
$\xi=\abs{\gw-\bar\gw}$ and $\eta=|\nabla'(\gw-\bar\gw)|$,
$$\sup\{\abs{D_1F(\gw,\nabla'\gw)}:(\gw,\nabla'\gw)\in G\}\leq\left\{\BA {ll} (p+q)\gg^{1-p}c_1^{p+q-1}\quad&\text{if }p\geq 1,
\\[2mm]
(p+q)\gg^{1-p}c_2^{p-1}c_1^{q}\quad&\text{if }0\leq p < 1.
\EA\right.
$$
and, since $|\gw||\nabla'\gw|\leq \frac{1}{2\gg}(\gg^2\gw^2+|\nabla'\gw|^2)$,
$$\sup\{\abs{D_2F(\gw,\nabla'\gw)}:(\gw,\nabla'\gw)\in G\}\leq\left\{\BA {ll} \myfrac{q}{2}\gg^{-p}c_1^{p+q-1}\quad&\text{if }p\geq 1,
\\[2mm]
\myfrac{q}{2}\gg^{-p}c_2^{p-1}c_1^{q}\quad&\text{if }0\leq p < 1.
\EA\right.
$$
\bel{R2'}\BA {lll}\myint{S^{n}}{}(\eta^2+\mu\xi^2)dS\leq
\gg^{1-p}c_*^{p+q-1}\myint{S^{n}}{}\left((p+q)\xi^2+\frac{q}{2}\gg\eta\xi\right)dS.
\EA\ee
Set
$$\Xi:=\left(\myint{S^{n}}{}\xi^2dS\right)^{\frac{1}{2}}\,\text{and }\; H:=\left(\myint{S^{n}}{}\eta^2dS\right)^{\frac{1}{2}},
$$
and recall that $c_*$ is defined in \eqref{cstar}.
Then $\sqrt{n\,}\Xi\leq H$. We define the polynomials
$$P(\Xi,H)=H^2-\frac{q}{2}\gg^{-p}c_*^{p+q-1}H\Xi+(\mu-(p+q)\gg^{1-p}c_*^{p+q-1})\Xi^2.
$$
and, putting $T=\frac{H}{\Xi}$  when $\Xi>0$, we have
\bel{R2''}\CP(T)=T^2-\frac{q}{2}\gg^{-p}c_*^{p+q-1}T+\mu-(p+q)\gg^{1-p}c_*^{p+q-1}.\ee
Then $T\geq \sqrt{n}$ since $n$ is the first nonzero eigenvalue of $-\Gd'$ in $H^1(S^n)$. \smallskip

\noindent Next we suppose that $\CP(\sqrt{n})\geq 0$. This means
\bel{R3}
n-\frac{q}{2}\gg^{-p}c_*^{p+q-1} \sqrt{n}+\mu-(p+q)\gg^{1-p}c_*^{p+q-1}\geq 0,
\ee
and it is equivalent to
\bel{R4}\BA {lll}
\left(\myfrac{q}{2}\gg^{-p}\sqrt{n}+(p+q)\gg^{1-p}\right)c_*^{p+q-1}\leq n+\mu.
\EA
\ee
We have three possibilities:\\

\noindent (i) either $\Xi>0$ and
\bel{R5}\BA {lll}
\myfrac{q^2}{4} \gg^{-2p}c_*^{2(p+q-1)}+4(p+q)\gg^{1-p}c_*^{p+q-1}\geq  4\mu,
\EA
\ee
then the polynomial $\CP$
admits two real roots $T_1\leq T_2$ and $\CP(T)\leq 0$. Jointly with the constraint on $T$ it means
$$T_1\leq T\leq T_2\leq \sqrt{n}\leq T.$$
Then $T= \sqrt{n}T$, which implies $\gw-\bar\gw=\gt\gf_1$ for some $\gt\in\BBR^*$. This is not compatible with the fact that $\gw$ solves \eqref{R0}. \\

\noindent (ii) either $\Xi>0$ and
\bel{R6}\BA {lll}
\myfrac{q^2}{4} \gg^{-2p}c_*^{2(p+q-1)}+(p+q)\gg^{1-p}c_1^{p+q-1}< 4\mu.
\EA
\ee
Then $\CP$ remains positive, which is impossible because of \eqref{R2'},  \\

\noindent (iii) or $\Xi=0$. In such a case $\gw=\bar\gw$, $\gw$ is a constant and $\nabla'\gw=0$. Therefore, if $(\ref{R4})$ holds $\gw=\bar\gw$ which ends the proof.\qeda\\

\noindent\Remark We notice that if we suppose $q=0$ in \eqref{R4} we find back condition (2.53) in \cite[Th 2.2]{LiVe}. %%%%%%%%%%%%%%%%%%%%%%%%%%%%%%%%%%%%%%%%%%%%%%%%%%%%%%%%%%%%%%%%%%%%BIFURCATION%%%%%%%%%%%%%%%%%%%%%%%%%%%%%%%%%%%%%%%%%%%%%%%%%%%%%%%%%%%%%%%%%%%%%%%%%%%%%%%%%%%%%%%%%%%%%%%%%%%%%%%%%%%%%%%%%%%%%%%%%%%%%%%%%%%%%%%%%%%%%%%%%%%%%%%%%%%%%%%%%%%%%%%%%%%%%%%%%%%%%%%%%%%%%%%%%%%%%%%%%%%%%%%%%%%%%%%
\subsection{Bifurcation}
In this section we are interested in solutions of $(\ref{L-1})$ which bifurcate from the constant solution $\gw_{\mu_*}$ defined by
\eqref{L-2} with $\mu_*=\frac{n}{p+q-1}$.
\bth{bifur} Assume $\gg>0$, $p+q-1>0$, and set $\mu_*=\frac{n}{p+q-1}$. Then there exists a neighborhood $\CO$ of
$(\mu_*,\gw_{\mu_*})$ in $\BBR\ti C^1(S^n)$ such that if $\gw$ is a solution of $(\ref{L-1})$ in $S^n$  satisfying $(\mu,\gw)\in \CO$, there holds either $(\mu,\gw)=(\mu,0)$ or $\mu>0$,  $(\mu,\gw)=(\mu_*+\ge(s),\gw_{\mu_*}+s(\phi_1+\phi(s)))$ where $s\mapsto \ge(s)$ is a $C^1$  function defined on $[0,\gt]$, vanishing at $s=0$ and $s\mapsto\phi(s)$ is a $C^1$ function defined on $[0,\gt]$, vanishing at $s=0$.
\es
\Proof
We set
$$\CL_\gm(\gw)=-\Gd'\gw+\mu\gw-\gw^p(\gg^2\gw^2+\abs{\nabla'\gw}^2)^{\frac q2}.
$$
We look for solutions under the form $\gw=\gw_{\mu_*}+\phi$ with $\gf$ small.
Then
$$\BA {ll}\gw^p(\gg^2\gw^2+\abs{\nabla'\gw}^2)^{\frac q2}=(\gw_{\mu_*}+\phi)^p\left(\gg^2(\gw_{\mu_*}+\phi)^2+\abs{\nabla'\gf}^2\right)^{\frac q2}\\[4mm]
\phantom{\gw^p(\gg^2\gw^2+)^{\frac q2}}
=\gg^q\gw_{\mu_*}^{p+q}(1+\myfrac{1}{\gw_{\mu_*}}\gf)^p\left(1+\myfrac{2}{\gw_{\mu_*}}\gf+\myfrac{1}{\gw_{\mu_*}^2}\gf^2+\myfrac{1}{\gg^2\gw_{\mu_*}^2}\abs{\nabla'\gf}^2\right)^{\frac q2}\\[4mm]
\phantom{\gw^p(\gg^2\gw^2+)^{\frac q2}}
=\gg^q\gw_{\mu_*}^{p+q}\left(1+\myfrac{p}{\gw_{\mu_*}}\gf+\myfrac{p(p-1)}{2\gw^2}\gf^2\right)\left(1+\myfrac{q}{\gw_{\mu_*}}\gf+
\myfrac{q(q-1)}{2\gw^2_*}\gf^2+\myfrac{q}{2\gg^2\gw_{\mu_*}^{2}}\abs{\nabla'\gf}^2\right)\\[4mm]
\phantom{\gw^p(\gg^2\gw^2+)^{\frac q2}---------}+O(\gf^3+\gf\abs{\nabla'\gf}^2)\\[4mm]
\phantom{\gw^p(\gg^2\gw^2+)^{\frac q2}}
=\gg^q\gw_{\mu_*}^{p+q}\left(1+\myfrac{p+q}{\gw_{\mu_*}}\gf+\myfrac{(p+q)(p+q-1)}{2\gw_{\mu_*}^2}\gf^2+\myfrac{q}{2\gg^2\gw_{\mu_*}^{2}}\abs{\nabla'\gf}^2\right)
\\[4mm]
\phantom{\gw^p(\gg^2\gw^2+)^{\frac q2}---------}+O(|\gf|^3+|\gf|\abs{\nabla'\gf}^2).
\EA$$
%%%%%%%%%%%%%%%%%%%%%%%%%%%%%%%%%%%%%%%%%%%%%%%%%%%%%%%%%%%%%%%%%%%%%%%%%%%%%%%%%%%%%%%%%%%%%%%%%%%%%%%%%%%%%%%%%%%%%%%%%%%%%
Since $\gg^{q}\gw_{\mu_*}^{p+q-1}=\mu_*$, we get
\bel{Bi-7}\BA {lll}\CL_\gm(\gw_{\mu_*}+\phi)= -\Gd'\gf+\mu\gw_{\mu_*}+\mu\gf
-\gg^q\gw_{\mu_*}^{p+q}\left(1+\myfrac{p+q}{\gw_{\mu_*}}\gf+\myfrac{(p+q)(p+q-1)}{2\gw_{\mu_*}^2}\gf^2+\myfrac{q}{2\gg^2\gw_{\mu_*}^{2}}\abs{\nabla'\gf}^2\right)\\[4mm]
\phantom{\CL(\gw_{\mu_*}+\phi)---------}+O(|\gf|^3+|\gf|\abs{\nabla'\gf}^2)
\\[4mm]
%-\ga\left(N-2+\ga-(p+q)\gg^{q-1}\gw_\mu^{p+q-1}\right)\gf
%
\phantom{\CL(\gw_\mu+\phi)}
=-\Gd'\gf+(\gm-\gm_*)\gw_{\mu_*}+\left(\gm-(p+q)\mu_*\right)\gf-\myfrac{(p+q)(p+q-1)}{2}\gg^q\gw_{\mu_*}^{p+q-2}\gf^2\\[4mm]
\phantom{\CL(\gw_\mu+\phi)---------}
-\myfrac{q}{2}\gg^{q-2}\gw_{\mu_*}^{p+q-2}\abs{\nabla'\gf}^2+O(|\gf|^3+|\gf|\abs{\nabla'\gf}^2).
\EA
\ee
Because $\mu_*(p+q-1)=n$, 
we can take $\gf=\ge\gf_1$ where $\ge$ is small and $\gf_1$ is the first non-zero eigenfunction (with corresponding eigenvalue $n$). Then
\bel{Bi-9-1}\BA {lll}
\CL_\gm(\gw_{\mu_*}+\ge\gf_1)= -\ge^2\gw_{\mu_*}^{p+q-2}\gg^{q-2}\left(\myfrac{(p+q)(p+q-1)}{2}\gg^2\gf_1^2
+\myfrac{q}{2}\abs{\nabla'\gf_1}^2\right)+O(\ge^3).
\EA
\ee
We want to apply \cite[Th 13.4, 13.5, Ex 2 p. 174]{Smo} and we consider solutions of
\bel{Bi-10-1}\BA {lll}
-\Gd'\gw+\mu\gw-|\gw|^{p-1}\gw(\gg^2\gw^2+|\nabla'\gw|^2)^{\frac{q}{2}}=0
\EA
\ee
depending only on the azimuthal angle $\gth_{n}:=\gth\in (0,\gp)$, which means
\bel{Bi-11*}\BA {lll}
-\gw_{\gth\gth}-(n-1)\cot_\gth\gw_\gth+\mu\gw-|\gw|^{p-1}\gw(\gg^2\gw^2+\gw_\gth^2)^{\frac{q}{2}}.
\EA
\ee
and denote by $C^{2,\gd}_{rad}(S^{n})$ ($\gd\in (0,1)$) the space of $C^{2,\gd}$ functions depending only on the angle $\gth$ (and thus radial with respect to the other variables $(\gth_1,...,\gth_{n-1})$). We recall that
$$\gw_{\mu_*}=(\gg^{-q}\mu_*)^{\frac{1}{p+q-1}}.
$$
Since the bifurcation point in \cite{Smo} are taken at $(\mu_*,\gw_{\mu_*})$, we put $\gw=\gw_{\mu_*}+w$ and
$$f(\mu,w)=-\Gd'w+\mu(\gw_{\mu_*}+w)-|\gw_{\mu_*}+w|^{p-1}(\gw_{\mu_*}+w)(\gg^2(\gw_{\mu_*}+w)^2+|\nabla' w|^2)^{\frac{q}{2}}.
$$
 Then
$$D_2f(\mu,0)=-\Gd'-(p+q-1)\mu I.
$$
If $(p+q-1)\mu=n\Longleftrightarrow\mu=\mu_*$, then ${\rm ker}D_2f(\mu_*,0)$ is spanned by $\gf_1:\gth\mapsto \cos\gth$ and
$${\rm R}(D_2f(\mu_*,0))=\left\{\psi\in C^{\gd}_{rad}(S^{n}):\myint{S^{n}}{}\psi\gf_1 dS=0\right\}.
$$
Finally  $D_1D_2f(\mu_*,0)(\mu,v)=-(p+q-1)\mu v$ thus $D_1D_2f(\mu_*,0)(\mu_* ,\gf_1)$ does not belong to ${\rm R}(D_2f(\mu_*,0))$. { Therefore the bifurcation theorem applies and there exists $\ge>0$ a $C^1$ curve $s\mapsto (\mu(s),\gf(s))$ defined on $[-\ge,\ge]$ with value in $\BBR\ti {\rm R}(D_2f(\mu_*,\gw_\mu))$ such that
\bel{Bi-12*}\BA {lcc}
(i)\qquad &\mu(0)=\mu_*&\qquad\qquad\\
(ii) \qquad &\gf(0)=0&\qquad\qquad\\
(iii)\qquad &f(\mu(s),s(\gf_1+\gf(s)))=0&\qquad\qquad
\EA\ee\smallskip

\noindent  Furthermore, there exists a neighborhood $\CO$ of $(\mu_*,\gw_{\mu_*})$ in which any solution of
$f(\mu,0)=0$ is either $(\mu,0)$ or under the form $(\mu(s),s(\gf_1+\gf(s)))$. Equivalently, any solution of
\eqref{Bi-11*} is either $(\mu,\gw_\mu)$ or is of the form $(\mu(s),\gw_\mu+s(\gf_1+\gf(s)))$}.

\noindent This is this last statement which applies in our case, which ends the proof. \medskip
\qeda

\noindent Remark. It would be interesting to study the direction of the bifurcation, which is not easy since the value $\gm^*$ is a second bifurcation, the first one occuring at $\gm=0$. We conjecture that the function $\ge(s)$ is positive. 

\mysection{Appendix}
%%%%%%%%%%%%%%%%%%%%%%%%%%%%%%%%%%%%%%%%%%%%%%%%%%%%%%%%%%%%%%%%%%%%%%%%%%%%%%%%%%%%%%
In this section we prove that under the assumptions of Theorem C we can choose the couple $(m,y)$ so that all our estimates in Section 3 are valid.

\subsection{Position of the problem}
We set for simplicity
\[
h=2e=(N-1)q\in\lbrack0,2(N-1)].
\]
The conditions to be satisfied by the parameters $(m,y)$ are:
\begin{equation}
A_{0}>0,\qquad B_{0}>0, \label{AB}%
\end{equation}
with $y\neq1$ and
with $y\neq1$ and
\begin{equation}
m<0,\qquad m+2+h>0 \label{two}
\end{equation}
and
\begin{equation}
(m+h+2)(2N-2-h)>(N-4-h)((N-1)p+h+1-N). \label{sigma}%
\end{equation}
We denote by $\mathcal{E}$ the ellipse of equation $\mathcal{E}(m,y)=0,$
where
\[
\mathcal{E}(m,y):=K{y}^{2}+2\,\left(  h+1\right)  (m-1)y-2\,\left(
h+1\right)  y+m(m-1)
\]
with
\[
K=\frac{2+h}{N}(N-1+Nh).
\]
Let $\mathcal{D}_{p}$ be the line defined by the equation $\mathcal{D}%
_{p}(m,y)=0,$ where
\[
\mathcal{D}_{p}(m,y)=y+am-bp
\]
with
\[
a=\frac{N+2+h}{2((N+1)h+N+2)},\qquad b=\frac{(N-1)(h+2)}{2((N+1)h+N+2)}%
\]
So the conditions $(\ref{AB})$ are equivalent to%

\[
\mathcal{E}(m,y)<0\text{ \quad and\quad\ }\mathcal{D}_{p}(m,y)>0,
\]
which means that the line $\mathcal{D}_p$ intersect the ellipse $\mathcal{E}$,
and $(m,y)$ lies inside $\mathcal{E}$ and above $\mathcal{D}_p$ where $p>0$.
First we write that $\mathcal{D}$ intersects $\mathcal{E},$ that means the
equation $\mathcal{E}(m,-am+bp)=0$ has at least one root. We obtain the equation $\mathcal{T}(m)=0$ with%
\[
\mathcal{T}(m)=\left(Ka^{2}-2\frac{bh}{N-1}\right)m^{2}-2\left(bp(Ka-1-h)-\frac{bh}%
{N-1}\right)m+bp(Kbp-2(1+h)).
\]
It is needed that its discriminant $\mathcal{J}$ be nonnegative and it is convenient to express it in terms of $t=bp$, hence
\[
\CJ=\left(tKa-\frac{bh}{N-1}-t(1+h)\right)^{2}-\left(Kt^{2}-2(1+h)t\right)\left(Ka^{2}-2\frac
{bh}{N-1}\right).
\]
We notice that
\begin{equation}
2a(1+h)-1=\frac{2bh}{N-1},\label{form}%
\end{equation}
therefore
$$\BA {lll}\displaystyle
\mathcal{J} =K^{2}t^{2}a^{2}+\left(\frac{bh}{N-1}+t(1+h)\right)^{2}%
-2tKa\left(\frac{bh}{N-1}+t(1+h)\right)-K^{2}t^{2}a^{2}\\[4mm]
\phantom{\mathcal{J}--------------- }
\displaystyle
+2(1+h)t\left(Ka^{2}-2\frac{bh}{N-1}\right)-2\frac{bh}%
{N-1}Kt^{2}\\[4mm]
\phantom{\mathcal{J} }
\displaystyle
=\left((1+h)^{2}-2aK(1+h)-2\frac{bh}{N-1}K\right)t^{2}+2b\left(  -\frac{h}%
{N-1}(1+h)-\frac{ahK}{N-1}\right.
\\[4mm]
\phantom{\mathcal{J} ----------------------}
\displaystyle
\left.+(1+h)K\frac{a^{2}}{b}\right)t  +\frac{b^{2}h^{2}%
}{(N-1)^{2}}.
\EA
$$
Hence $\mathcal{J}=-\frac{b^{2}}{N}\tilde{G}(p,h)$ where
\bel{fun}\BA {lll}\displaystyle
\tilde{G}(p,h)=(K-(h+1)^{2})p^{2}+2\left(\frac{h(1+h)}{N-1}%
-\frac{a^{2}}{2b}\right)  p+\frac{h^{2}}{(N-1)^{2}}\\[4mm]\displaystyle
\phantom{\tilde{G}(p,h)}
=(N-1)((N-1)h+N-2)p^{2}\\[4mm]\displaystyle
\phantom{\tilde{G}(p,h)--------}+\left(  N{h}^{2}-({N}^{2}+N-1)h-{N}^{2}-N+2\right)
p-\frac{N}{N-1}h^{2}.%
\EA
\ee
So we find precisely that $\tilde{G}(p,(N-1)q)=G(p,q)$ where the function $G$
is given in Theorem C. The equation $\tilde{G}(p,h)=0$ has two
roots with opposite sign in $p,$ that we call $=p_{0}(q)>0\geq p_{1}(q)$. Both correspond
to the fact that the  lines $\CD_{p_i(q)}$ are tangent to the ellipse $\mathcal{E}$.
The region $\CD_{p_{1}(q)}(m,y)>0$ contains the whole region $\CE(m,y)<0$ (the interior of $\CE$) while $\CD_{p_{0}(q)}(m,y)>0$
 has an empty intersection with the region $\CE(m,y)<0$. Hence for $p>0,$ the line $\mathcal{D}_{p}$
intersect the ellipse if and only if $p<p_{0}(q)$. If $q=0,$ then we find that
$p_{0}(0)=\frac{N+2}{N-2},$ which was the precise optimal value obtained for
the Emden-Fowler equation -$\Delta u=u^{p}$. Now for $p=p_{0}(q),$ the line
$\mathcal{D}_{p_{0}(q)}$ is tangent to the ellipse at some point $(m_{0}%
,y_{0})=(m_{0}(q),y_{0}(q))$ in the upper part of $\mathcal{E},$ given by
\[
y_{0}=\frac{1}{K}\left(  (1-m_{0})(h+1)+\sqrt{(1-m_{0})((1-m_{0}%
)(h+1)^{2}+m_{0}K}\right)  >0.
\]

Suppose that we have proved that $m_{0}(q)$ satisfies the conditions
$(\ref{two})$ and $(\ref{sigma})$, then for a given $p<p_{0}(q)$ any couple $(y,m)$
with $0<y<y_{0}(q),$ $y\neq1,$ and $m=m_{0}(q)$ will satisfy all the required
conditions. Therefore it is sufficient to prove $(\ref{two})$ and $(\ref{sigma})$ in case
$m=m_{0}(q).$\medskip

\nind \Remark
If it happens that $m_{0}+2+h=0,$ then we can take $m=m_{0}%
+\ge$ with $\ge>0$ small enough such that $(m,y)$ stays in
$\mathcal{E}.$ We know from \cite{BV-V} that it
happens precisely when $N=3$ and $h=0.$ We will see below that it is the only case.\medskip

Next we compute $m_{0}.$ We note here that the discriminant of $p\mapsto \tilde G(p,h)$ is
\[
\mathcal{H}=\left(  N{h}^{2}-({N}^{2}+N-1)h-{N}^{2}-N+2\right)  ^{2}%
+4Nh^{2}((N-1)h+N-2)>0
\]
and it can be written under the form $\mathcal{H}=(Nh+N-1)\mathcal{M}$, where
\[
\mathcal{M=M}(h)=N{h}^{3}-\left(  2\,{N}^{2}-N+1\right)  {h}^{2}+\left(
{N}^{3}+2\,{N}^{2}-2\,N-4\right)  h+{(N-1)(N+2)}^{2}.%
\]
Then $p_{0}=p_{0}(q)$ is given by
\begin{equation}\label{defip}
p_{0}=\frac{-\left(  N{h}^{2}-({N}^{2}+N-1)h-{N}^{2}-N+2\right)
+\sqrt{(Nh+N-1)\mathcal{M}}}{2(N-1)((N-1)h+N-2)}.%
\end{equation}
Since $m_{0}$  is also the minimizer of the trinomial $m\mapsto \mathcal{T}(m)$, it is expressed by
\[
m_{0}=b\frac{(N-1)(Ka-1-h)p_{0}-h}{(N-1)Ka^{2}-2bh},
\]
and we obtain, after some computation,
\[
m_{0}=\frac{Q_{1}(h)+(N-1)Q_{2}(h)p_{0}}{\mathcal{M}},
\]
where%
\begin{equation}\label{Q1}
Q_{1}(h)=-2\,Nh\left(  \left(  N+1\right)  h+N+2\right),%
\end{equation}%
and
\begin{equation}\label{Q2}
Q_{2}(h)=N{h}^{3}-\left(  {N}^{2}-3\,N+1\right)  {h}^{2}-\left(  {N}%
^{2}-N+4\right)  h-2\,(N-2).%
\end{equation}
Replacing $p_{0}$ by its value given in (\ref{defip}), we deduce after some
simplifications that
\[
2((N-1)h+N-2)m_{0}\mathcal{M}=-(Nh^{2}+(N+1)h+2)\mathcal{M}+Q_{2}%
\sqrt{(Nh+N-1)\mathcal{M}}.
\]
Hence%
\begin{equation}
2((N-1)h+N-2)m_{0}=-P_{1}+Q_{2}\sqrt{\frac{Nh+N-1}{\mathcal{M}}}\label{valm}%
\end{equation}
with
\[
P_{1}(h)=Nh^{2}+(N+1)h+2.
\]
Note that $m_{0}$ can be also obtained equivalently by expressing the fact that
$(m_{0},y_{0})$ belongs in the upper part of $\mathcal{E}$ and the slope of its tangent here has value $-a.$\medskip

\nind\Remark
When $q=0,$ we rediscover the values given in \cite{BV-V},
\begin{equation}
p_{0}=\frac{N+2}{N-2},\quad m_{0}=-\frac{2}{N-2},\quad y_{0}=\frac{N}%
{N-2}.\label{qzero}%
\end{equation}

%%%%%%%%%%%%%%%%%%%%%%%%%%%%%%%%%%%%%%%%%%%%%%%%%%%%%%%%%%%%%%%%%%%%%%%%%%%%%%%%%%%%%%%%%%%%%%%%%%%%%%%%%%%%%%%%%%%%%%%%%%%%%%%%%%%%%%%%%%%%%%%%%%%%%%%%%%%%%%%%%%%%%%%%%%%
%%%%%%%%%%%%%%%%%%%%%%%%%%%%%%%%%%%%%%%%%%%%%%%%%%%%%%%%%%%%%%%%%%%%%%%%%%%%%%%%%%%%%%%%%%%%%%%%%%%%%%%

\subsection{Proof that $m_0<0$ for $h\in [0,2(N-1)]$}

We present a proof which avoids the heavy computation of $m_{0}$. The point
$(m_{0},y_{0})$ belongs to the upper part of $\mathcal{E},$ where by concavity the
slope $m\mapsto y^{\prime}(m)$ is a decreasing function. Hence, the claim will
follow provided $y^{\prime}(m_{0})=-a>y^{\prime}(0)$. We obtain
directly $y(0)=\frac{2(1+h)}{K}$ and
\[
y^{\prime}(0)\left[  2Ky(0)-2(1+h)\right]  )+2(1+h)y(0)-1=0,
\]%
\[
\left\vert y^{\prime}(0)\right\vert =\frac{1}{2(1+h)}(\frac{4(1+h)^{2}}%
{K}-1);
\]
hence, from $(\ref{form})$,%
\begin{align*}
(1+h)(\left\vert y^{\prime}(0)\right\vert -a)  &  =\frac{2(1+h)^{2}}{K}%
-\frac{1}{2}-a(1+h)=\frac{2(1+h)^{2}}{K}-(1-\frac{bh}{N-1})\\
&  =\frac{2(1+h)^{2}}{K}-\frac{N+2+(N-1)h-h^{2}}{N+2+(N+1)h}.%
\end{align*}
It is therefore required that
\begin{align*}
&  2N(1+h)^{2}(N+2+(N+1)h)-(2+h)(Nh+N-1)(N+2+(N-1)h-h^{2})\\
&  =Nh^{4}+(N^{2}+6N-1)h^{3}+(2N^{2}+12N-3)h^{2}+(N^{2}+9N)h+2(N+2)>0
\end{align*}
which clearly holds.
%%%%%%%%%%%%%%%%%%%%%%%%%%%%%%%%%%%%%%%%%%%%%%%%%%%%%%%%%%%%%%%%%%%%%%RIGIDITY%%%%%%%%%%%%%%%%%%%%%%%%%%%%%%%%%%%%%%%%%%%%%%%%%%%%%%%%%%%%%%
%%%%%%%%%%%%%%%%%%%%%%%%%%%%%%%%%%%%%%%%%%%%%%%%%%%%%%%%%%%%%%%%%%%%%%%%%%%%%%%%%%%%%%%%%%%%%%%%%%%%%%%

\subsection{Proof that $m_0+2+h>0$ for $h\in (0,2(N-1)]$}
From $(\ref{valm})$ the value of $m_{0}+2+h$ is given by
\begin{equation}
2((N-1)h+N-2)(m_{0}+2+h)=P_{2}(h)+Q_{2}(h)\sqrt{\frac{Nh+N-1}{\mathcal{M}}},
\label{cro}%
\end{equation}%
with
\begin{equation}
P_{2}(h)=((N-2){h}^{2}+(5\,N-9)h+4\,N-10, \label{P1}%
\end{equation}
and where we recall that
\[
Q_{2}(h)=N{h}^{3}-\left(  {N}^{2}-3\,N+1\right)  {h}^{2}-\left(  {N}%
^{2}-N+4\right)  h-2\,N-4.
\]
Note that $P_{2}(h)>0$ for any $h\in\lbrack0,2N-2]$. Then $m_{0}+2+h>0$ as
soon as $Q_{2}(h)>0$. Since $Q_2$ can be written under the form
\[
Q_{2}(h)=(h+2)(Nh-N-1)(h-N+2)-2N^{2},%
\]
it is an increasing function of $h$ on $[N-2,2N-2]$.\smallskip

\nind {\it The case }$N\geq4$. Here $Q_{2}(N-1)>0$, thus
$Q_{2}(h)>0$ on $[N-1,2N-2]$. Therefore it is sufficient to prove the assertion when
$h\in\lbrack0,N-1]$. Our aim is to prove that there holds
\[
P_{2}\sqrt{R}+Q_{2}>0\text{ where }R=\frac{\mathcal{M}}{Nh+N-1}%
\]
in this interval. By division we obtain, since $Nh+N-1<N^{2}-1,$
\bel{fre}\BA {lll}
R =(N-h)^{2}+6N+2-\myfrac{2(N+2)(N+1)}{N}+\myfrac{4(N-1)(2N+1)}{N(Nh+N-1)}%
\\[4mm]\displaystyle\phantom{R}
\geq(N-h)^{2}+6N+2-\frac{2(N+2)(N+1)}{N}+\frac{4(2N+1)}{N(N+1)}
\\[4mm]\displaystyle\phantom{R}=(N-h)^{2}+4(N-1)+\frac{4}{N+1}.
\EA\ee
In particular $\sqrt{R}\geq N-h,$ thus $P_{1}\sqrt{R}+Q_{2}\geq S$
with
\[
S(h)=(N-h)P_{2}+Q_{2}=2\,{h}^{3}-4\left(  N-2\right)  {h}^{2}+\left(
4\,{N}^{2}-12\,N+6\right)  h+4\,({N}^{2}-3\,N-1).
\]
We see that $S(0)>0$ because $N\geq4,$ and since
\[
S^{'}(h)=6h^{2}-2(4N-8)h+4N^{2}-12N+6>0
\]
is positive for any $h\in\BBR$, $S$ is an increasing function of $h$. This yields $S(h)>0$ on $h\in\lbrack0,N-1]$ and completes the proof in this case. \smallskip

\nind{\it The case }$N=3$. Here we cannot use the minorization of $R$
since equality holds for $h=0$, as it was noticed above. Also, we observe
that the cubic polynomial $Q_{2}(h)=3h^{3}-h^{2}-10h-10$ has its largest root $h_{0}$ in the
interval $(2,3)$, since $Q_{2}(2)=-10$ and $Q_{2}(3)=32$ and $Q_2'(h)>0$ on $[2,\infty)$. Hence we need only to
prove the inequality when $h\in(0,3).$ For this aim, it is sufficient
that $P_{2}\sqrt{R}+Q_{2}>0$, which will be ensured provided $\mathcal{M}%
P_{2}^{2}-(3h+2)Q_{2}^{2}>0$. After some computation it reduces to prove that
\[
h\left(  -6\,{h}^{6}+5\,{h}^{5}+38\,{h}^{4}+35\,{h}^{3}+337\,{h}%
^{2}+484\,h+160\right)  >0.
\]
This inequality is clearly true, since $6\,{h}^{4}-38\,{h}^{2}-337<0$ on $(0,3)$. So
finally $m_{0}+2+h>0$ for any $h\in\left(  0,4\right]  .$

%%%%%%%%%%%%%%%%%%%%%%%%%%%%%%%%%%%%%%%%%%%%%%%%%%%%%%%%%%%%%%%%%%%%%%%%%%%%%%%%%%%%%%%%%%%%%%%%%%%%%%%%%%%%%%%%%%%%%%%%%%%%%%%%%%%%%%%%%%%%%%%%%%%%%%%%%%%%%%%%%%%%%%%%%%%
%%%%%%%%%%%%%%%%%%%%%%%%%%%%%%%%%%%%%%%%%%%%%%%%%%%%%%%%%%%%%%%%%%%%%%%%%%%%%%%%%%%%%%%%%%%%%%%%%%%%%%%
\subsection{Proof that $\gs>\frac N2$ for $N\geq 3$}
We have to prove that for any $h\in\left(  0,2(N-1)\right]$ there holds,%
\bel{Ap1}
(m_{0}+h+2)(2N-2-h)>(N-4-h)((N-1)(p_{0}-1)+h).
\ee
In the preceding
step we have already shown that the left hand side is positive and that $(N-1)(p_{0}-1)+h)=(N-1)(p+q-1)>0$, inequality $(\ref{Ap1})$ is valid
for $N=3,4$ or $N\geq5$ and $h\in\left[  N-4,2(N-1)\right]$. Henceforth we assume $N\geq5$ and $h\in\left(  0,N-4\right)$.
By replacing the value of $m_{0}(h)+h+2$ given the preceding section and $p_{0}$ by its value,
given in $(\ref{defip})$, we obtain after some computation that the relation can be
expressed under the form
\[
Q_{3}\sqrt{R}+Q_{4}>0,
\]
where
$$\BA {lll}
Q_{3}(h)=\left(  2\,N-2\right)  {h}^{2}+\left(  {N}^{3}-2\,{N}^{2}%
+10\,N-6\right)  h+{N}^{3}-3\,{N}^{2}+6\,N-4),
\EA$$%
and
$$\BA {lll}
Q_{4}(h)=\left(  {N}^{3}+2\,{N}^{2}-6\,N-2\right)  {h}^{2}%
\\[1mm]
\phantom{------}-\left(  {N}^{4}-{N}^{3}-17\,{N}^{2}+12\,N+8\right)  h-(N-1)({N}^{3}-8N-8).
\EA$$
%%%%%%%
From $(\ref{fre})$, we get
$$\BA {lll}\displaystyle
R\geq(N-h)^{2}+4(N-1)=(N-h)^{2}\left(  1+\frac{4(N-1)}{(N-h)^{2}}\right)\\[4mm]\displaystyle
\phantom{R\geq(N-h)^{2}+4(N-1)}
\geq(N-h)^{2}\left(  1+\frac{4(N-1)}{N^{2}}\right).
\EA$$
Setting $\tau:=\sqrt{1+\frac{4(N-1)}{N^{2}}}$ , it is therefore sufficient that
$$
Q_{5}:=\tau(N-h)Q_{3}+Q_{4}>0.
$$
An explicit computation yields
$$\BA {lll}
Q_{5}(h)=-2\left(  \,N-1\right)  \tau{h}^{3}+\left[\left(  1-\tau\right)  {N}%
^{3}+\left(  4\,\tau+2\right)  {N}^{2}-6\left(  2\,\tau+1\right)
N+6\,\tau-2\right]{h}^{2}
\\[1mm]\phantom{Q_{5}(h)}
+\left[  \left(  \tau-1\right)  {N}^{4}+\left(  1-3\,\tau\right)  {N}%
^{3}+\left(  13\,\tau+17\right)  {N}^{2}-12\left(  \,\tau+1\right)
N+4\,(\tau-2)\right]  h
\\[1mm]\phantom{Q_{5}(h)}
+\left(  \tau-1\right)  {N}^{4}+\left(  1-3\tau\right)  {N}^{3}+\left(
6\,\tau+8\right)  {N}^{2}-4\,\tau N-8
\EA$$
The function $Q_{5}$ is a cubic with a negative leading coefficient. We claim that it is
positive for $h=0$ and for $h=N-4$ and increasing near $0$. Indeed%
$$\BA {lll}
Q_{5}(0)  =\left(  \tau-1\right)  {N}^{4}+\left(  1-3\tau\right)  {N}^{3}+\left(
6\,\tau+8\right)  {N}^{2}-4\,\tau N-8\\[1mm]\phantom{Q_{5}(0)  }
 =(N-1)(\tau N({N}^{2}-2N+4)-(N^{3}-8N-8)).
\EA$$
Since $\tau^2=\frac{N^2+4N-4}{N^2}$, the inequality $Q_{5}(0)>0$ is equivalent to
$$
(N^{2}+4N-4)({N}^{2}-2N+4)^{2}-(N^{3}-8N-8)^{2}>0,
$$
which can be easily verified since $N\geq5$. Next
$$\BA {lll}
Q_{5}(N-4)  =\tau(4N^{4}-12N^{3}-12N^{2}+32N-48)-2N^{4}+8N^{3}%
+6N^{2}-40N-8\\[1mm]
\phantom{Q_{5}(N-4)}
>2N^{4}-4N^{3}-6N^{2}-8N-56>0
\EA$$
For the derivative, we compute
\begin{align*}
Q_{5}^{\prime}(0) &  =\left(  \tau-1\right)  {N}^{4}+\left(  1-3\,\tau\right)
{N}^{3}+\left(  13\,\tau+17\right)  {N}^{2}-12\left(  \,\tau+1\right)
N+4\,(\tau-2)\\
&  =\tau(N^{4}-3{N}^{3}+13N^{2}-12N+4)-({N}^{4}-N^{3}-17N^{2}+12N+8).
\end{align*}
Replacing $\gt$ by its value, the sign of $Q_{5}^{\prime}(0)$ is the same as the one of
$$\BA {lll}
(N^2+4N-4)\left(N^{4}-3{N}^{3}+13N^{2}-12N+4\right)^2-N^2\left({N}^{4}-N^{3}-17N^{2}+12N+8\right)^2
\EA$$
which is equal to
\begin{align*}
40N^{8}+4N^{7}-580N^{6}+1492N^{5}-1964N^{4}+2432\allowbreak N^{3}%
-1424N^{2}+448N-64
\end{align*}
and is clearly positive since $N\geq 5$. At end $Q_{5}$ stays positive on $\left[  0,N-4\right]$ and the proof is achieved.

%%%%%%%%%%%%%%%%%%%%%%%%%%%%%%%%%%%%%%%%%%%%%%%%%%%%%%%%%%%%%%%%%%%%%%%%%%%%%%%%%%%%%%%%%%%%%%%%%%%%%%%
%%%%%%%%%%%%%%%%%%%%%%%%%%%%%%%%%%%%%%%%%%%%%%%%%%%%%%%%%%%%%%%%%%%%%%%%%%%%%%%%%%%%%%%%%%%%%%%%%%%%%%%
\subsection{Comparison of the regions of Theorems B, D and Theorem C}
In the variable $h:=(N-1)q$ the  curves for Theorems B, D are given by
\begin{align}
(N-1)p+h &  =N+3\qquad\text{for\quad}p>1,\label{pre}\\
(N-1)p+h &  =N-1+{\frac{\left(  p+1\right)  ^{2}}{p}}\qquad\text{for\quad
}p\in [0,1].\label{seco}
\end{align}
In order to show that the curve $\tilde{G}(p,h)=0,$ where $\tilde{G}$ is
defined at $(\ref{fun})$, is above them, we only need to show that $\tilde
{G}(p,h)<0$ on these curves. For the first curve defined by $(\ref{pre})$, using
that $(h-2)(h-3)\geq-\frac{1}{4}$ for all $h,$ we obtain that for
$h\in\left[  0,2(N-1)\right],$
\begin{align*}
(N-1)G\left(\frac{N+3-h}{N-1},h\right) &  =-{h}^{3}+3\,{h}^{2}%
-4\,N+4\,h-12\\&
=-(h+2)(h-2)(h-3)-4N\\
&  \leq\frac{h+2}{4}-4N\\&\leq\frac{N}{2}-4N<0.
\end{align*}
As for the second curve, we check that
$$\BA {lll}\displaystyle
(N-1)\frac{p^{2}}{(p+1)^{2}}\tilde{G}\left(p,N-1+{\frac{\left(  p+1\right)  ^{2}%
}{p}}-(N-1)p\right)\\[4mm]\phantom{---------}
=-p(p-1)^{2}N^{2}+(3p^{3}-2p^{2}-p-1)N-p^{2}(2p+1)<0
\EA
$$
for $p\in\lbrack0,1],$ because $3p^{3}-2p^{2}-p-1)N-p^{2}(2p+1)\leq
p^{2}-p-1\leq-1.$

\subsection{Final remark about the parameter $\gb$}
When we set $u=v^{-\beta}$ then $y=\frac{\beta+1}{\beta}$ a natural question
 is about the sign of $\beta$. We have seen that $y_{0}>0,$ and when $p<p_{0}$
we have chosen $y=y_{0}$ if $y_{0}\neq1,$ or $y=y_{0}-\ge$ in case
$y_{0}=1$ (or in the special case $N=3,$ where $m_{0}=-2,$ so we have taken
$m=-2+\ge).$ So either $y_{0}>1,$  $\beta>0;$ or $y_{0}\leq1,$
$\beta<0.$ We can remark that for $q=0$, we have $\beta>0,$ since
$y_{0}=\frac{N}{N-2}$ from (\ref{qzero}). But for $q=2$ we find
\[
p_{0}=\frac{4}{2N-3}=-m_{0},\qquad y_{0}=\frac{2}{2N-3}%
\]
hence $y_{0}<1,$ thus $\beta<0.$

%%%%%%%%%%%%%%%%%%%SYMMETRY%%%%%%%%%%%%%%%%%%%%%%%%%%%%%%%%%%%%%%%%%%%%%%%%%%%%%%%%%%%%%%%%%%%%%%%%%%%%%%%%%%%%%%%%%%%%%%%%%%%%%%%%%%%%%%%%%%%%%%%%%%%%%%%%%%%%%%%%%%%%%%%%%%

%%%%%%%%%%%%%%%%%%%%%%%%%%%%%%%%%%%%%%%%%%%%%%%%%%%%%%%%%%%%%%%%%%%%%%%%%%%%%Rigidity on $S^{N-1}_+$%%%%%%%%%%%%%%%%%%%%%%%%%%
%%%%%%%%%%%%%%%%%%%%%%%%%%%%%%%%%%%%%%%%%%%%%%%%%%%%%%%%%%

%%%%%%%%%%%%%%%%%%%%%%%%%%%%%%%%%%%%%%%%%%%%%%%%%%%%%%%%%%%%%%%%%%%%%%%%%%%%%%%%%%%%%%%%%%%%%%%%%%%%%%%%%%%%%%%%%%%%%%%%%%%%%%%%%%%%%%%%%%%%%%%%%%%%%%%%%%%%%%%%%%%%%%%%%%%%%%%%%%%%%%%%%%%%%%%%%%%%%%%%%%%%%%%%%%%%%%%%
%%%%%%%%%%%%%%%%%%%%%%%%%%%%%%%%%%%%%%%%%%%%%%%%%%%%%%%%%%%%%%%%%%%already obtained by Borghol and V\'eron \cite{BoVe1} in the case $p=N$%%%%%%%%%%%%%%%%%%%%%%%%%%%%%%%%%%%%%%%%%%%%%%%%%%%%%%%%%%%%%%%%%%%%%%%%%%%%%%%%%%%%%%%%%%%%%%%%%%%%%%%%%%%%%%%%%%%%%%%%%%
%%%%%%%%%%%%%%%%%%%%%%%%%%%%%%%%%%%%%%%%%%%%%%%%%%%%%%%%%%%%%%%%%%%%%%%%%%%%%%%%%%%%%%%%%%%%%%%%%%%%%%%%%%%%%%%%%%%%%%%%%%%%%%%%%%%%%%%%%%%%%%%%%%%%%%%%%%%%%%%%%%%%%%%%%%%%%%%%%%%%%%%%%%%%%
%%%%%%%%%%%%%%%%%%%%%%%%%%%%%%%%%%%%%%%%%%%%%%%%%%%%%%%%%%%%%%%%%%%%%%%%%%%%%%%%%%%%%%%%%%%%%%%%%%%%%%%%%%%%%%%%%%%%%%%%%%%%%%
%%%%%%%%%%%%%%%%%%%%%%%%%%%%%%%%%%%%%%%%%%%%%%%%%%%%%%%%%%%%%%%%%%%%%%%%%%%%%%%%%%%%%%%%%%%%%%%%%%%%%%%%%%%%%%%%%%%%%%%%%%%%%%

\end{document}